\documentclass[12pt,a4paper]{article}
\usepackage{amssymb,amsmath,amsfonts,amscd}

\newcommand{\initpar}{\makebox[\parindent]{ }}
\newtheorem{twX}{Theorem}[section]
\newtheorem{lemX}[twX]{Lemma}
\newtheorem{propX}[twX]{Proposition}
\newtheorem{defX}[twX]{Definition}
\newtheorem{uwX}[twX]{Remark}
\newtheorem{wnX}[twX]{Corollary}
\newtheorem{exaX}[twX]{Example}
\newtheorem{pytX}[twX]{Question}

\def\kropka{\hspace{-2mm}\mbox{\bf .}}

\def\kon{\mbox{$\square$}}

\newcommand{\btw}[1]{\begin{twX}\label{#1}\kropka}
\newcommand{\bbtw}[2]{\begin{twX}[#1]\label{#2}\kropka\hspace{2mm}}
\newcommand{\etw}{\end{twX}}

\newcommand{\blem}[1]{\begin{lemX}\label{#1}\kropka}
\newcommand{\bblem}[2]{\begin{lemX}[#1]\label{#2}\kropka\hspace{2mm}}
\newcommand{\elem}{\end{lemX}}

\newcommand{\bprop}[1]{\begin{propX}\label{#1}\kropka}
\newcommand{\bbprop}[2]{\begin{propX}[#1]\label{#2}\kropka\hspace{2mm}}
\newcommand{\eprop}{\end{propX}}

\newcommand{\bdefn}[1]{\begin{defX}\label{#1}\kropka}
\newcommand{\bbdefn}[2]{\begin{defX}[#1]\label{#2}\kropka\hspace{2mm}}
\newcommand{\edefn}{\end{defX}}

\newcommand{\buw}[1]{\begin{uwX}\label{#1}\kropka}
\newcommand{\bbuw}[2]{\begin{uwX}[#1]\label{#2}\kropka\hspace{2mm}}
\newcommand{\euw}{\end{uwX}}

\newcommand{\bwn}[1]{\begin{wnX}\label{#1}\kropka}
\newcommand{\bbwn}[2]{\begin{wnX}[#1]\label{#2}\kropka\hspace{2mm}}
\newcommand{\ewn}{\end{wnX}}

\newcommand{\bexa}[1]{\begin{exaX}\label{#1}\kropka}
\newcommand{\bbexa}[2]{\begin{exaX}[#1]\label{#2}\kropka\hspace{2mm}}
\newcommand{\eexa}{\end{exaX}}

\newcommand{\bpyt}[1]{\begin{pytX}\label{#1}\kropka}
\newcommand{\bbpyt}[2]{\begin{pytX}[#1]\label{#2}\kropka\hspace{2mm}}
\newcommand{\epyt}{\end{pytX}}

\newenvironment{dowod}
 {\begin{bf}Proof. \end{bf}}{\kon \medskip \par}

\newcommand{\f}{\varphi}

\newcommand{\MOD}[3]{#1\equiv#2\hspace{-1mm}\pmod{#3}}

\newcommand{\E}{\varepsilon}

\begin{document}
\thispagestyle{empty}

\begin{center}
\begin{LARGE}

{\bf Constants of cyclotomic derivations}

\end{LARGE}

\bigskip

{\bf Jean Moulin  Ollagnier$^1$} and {\bf Andrzej Nowicki$^2$}

\bigskip

{$^1$\small
            Laboratoire LIX, \'Ecole Polytechnique,
            F 91128 Palaiseau Cedex, France,\\
            (e-mail :  Jean.Moulin-Ollagnier@polytechnique.edu).}

\medskip

{\small $^2$Nicolaus Copernicus University, Faculty of Mathematics and
Computer Science,\\
 87-100 Toru\'n, Poland, (e-mail:  anow@mat.uni.torun.pl).}

\bigskip


\end{center}

\medskip

\begin{abstract}
Let $k[X]=k[x_0,\dots,x_{n-1}]$ and $k[Y]=k[y_0,\dots,y_{n-1}]$
be the polynomial rings in $n\geqslant 3$ variables
over a field $k$ of characteristic zero
containing  the $n$-th roots
of unity. Let $d$ be the cyclotomic derivation of $k[X]$, and let $\Delta$
be the factorisable derivation of $k[Y]$
associated with $d$, that is, $d(x_j)=x_{j+1}$ and $\Delta(y_j)=y_j(y_{j+1}-y_j)$ for all $j\in\mathbb Z_n$.
We describe polynomial constants and rational constants of these derivations.
We prove, among others, that the field of constants of $d$ is
a field of rational functions over $k$ in $n-\f(n)$ variables, and that the ring of constants
of $d$ is a polynomial ring if and only if $n$ is a power of a prime.
Moreover, we show that the ring of constants of $\Delta$ is always
equal to $k[v]$, where $v$ is the product $y_0\cdots y_{n-1}$,
and we describe the field of constants of $\Delta$
in two cases: when $n$ is power of a prime, and when $n=p q$.
\footnotetext{
Corresponding author : Andrzej Nowicki, Nicolaus Copernicus
University, Faculty of Mathematics and Computer Science, ul. Chopina 12/18,
87--100 Toru\'n, Poland.\quad E-mail:  {\tt anow@mat.uni.torun.pl}.}
\end{abstract}

\medskip


\medskip


\medskip

\section*{Introduction}
\initpar
Throughout this paper $n\geqslant3$ is an integer,
$k$  is a field of characteristic zero containing  the $n$-th roots
of unity,
and $k[X]=k[x_0,\dots,x_{n-1}]$ and
$k[Y]=k[y_0,\dots,y_{n-1}]$ are polynomial rings over $k$
in $n$ variables.
We denote by $k(X)=k(x_0,\dots,x_{n-1})$ and $k(Y)=k(y_0,\dots,y_{n-1})$
the fields of quotients of $k[X]$ and $k[Y]$, respectively.
We fix the notations $d$ and $\Delta$ for the following two derivations, which
we call {\it cyclotomic derivations}.
We denote by $d$ the derivation of $k[X]$ defined by
$$
d(x_j)=x_{j+1},\quad\mbox{for}\quad j\in\mathbb Z_n,
$$
and we denote by $\Delta$ the  derivation of $k[Y]$ defined by
$$
\Delta(y_j)=y_j(y_{j+1}-y_j),\quad\mbox{for}\quad j\in\mathbb Z_n.
$$
We denote also by $d$ and $\Delta$ the unique extension of~$d$ to $k(X)$
and the unique extension of $\Delta$ to $k(Y)$, respectively.
We will show that there are some important relations between $d$ and $\Delta$.
In this paper we study polynomial and rational
constants of these derivations.

\smallskip

In general, if $\delta$ is a derivation of a commutative $k$-algebra $A$,
then we denote by $A^\delta$ the $k$-algebra of constants of $\delta$, that is,
$
A^\delta=\left\{ a\in A; \ \delta(a)=0\right\}.
$
For a given derivation $\delta$ of $k[X]$,
we are interested in some descriptions
of $k[X]^\delta$ and $k(X)^\delta$.
However, we know that such descriptions are usually
difficult to obtain.
Rings and fields of constants appear in various classical problems;
for details we refer to \cite{Essen00}, \cite{Freu}, \cite{NN}  and \cite{Now}.
The mentioned problems are already difficult for factorisable derivations.
We say that a derivation $\delta:k[X]\to k[X]$ is {\it factorisable} if
$$\delta(x_i)=x_i\sum_{j=0}^{n-1} a_{ij} x_j$$
for all $i\in \mathbb Z_n$, where each $a_{ij}$ belongs to $k$.
Such factorisable derivations and factorisable systems of ordinary
differential equations were intensively studied from a long time;
see for example
\cite{HofbauerSigmund88}, \cite{GrammaticosMoulin90}, \cite{MoulinNS} and \cite{NowUJ}.
Our derivation $\Delta$ is factorisable, and the derivation $d$ is {\it monomial},
that is, all the polynomials $d(x_0),\dots,d(x_{n-1})$ are monomials.
With any given monomial derivation $\delta$ of $k[X]$
we may associate, using a special procedure,
the unique factorisable derivation~$D$ of $k[Y]$ (see \cite{Mac}, \cite{NowZ}, \cite{MoulinN},
for details), and then,
very often, the problem of descriptions of
$k[X]^\delta$ or $k(X)^\delta$ reduces to the same problem for the
factorisable derivation~$D$.

\smallskip

Consider a derivation $\delta$ of $k[X]$ given by $\delta(x_j)=x_{j+1}^s$ for $j\in\mathbb Z_n$,
where $s$ is an integer.
Such $d$ is called a Jouanolou derivation (\cite{Jou}, \cite{MoulinNS}, \cite{Mac}, \cite{Zol}).
The factorisable derivation $D$, associated with this $\delta$,
is a derivation of $k[Y]$
defined by
$D(y_j)=y_j(sy_{j+1}-y_j),$
for $j\in\mathbb Z_n$.
We proved in \cite{Mac} that if $s\geqslant2$ and $n\geqslant3$ is prime,
then the field of constants of $\delta$ is trivial, that is, $k(X)^\delta=k$.
In $2003$ H. \.Zo\l \c{a}dek \cite{Zol}
proved the for $s\geqslant2$, it is also true
for arbitrary $n\geqslant3$;
without the assumption that $n$ is prime.
The central role, in his and our proofs, played
some extra properties of the associated  derivation~$D$.
Indeed, for $s\geqslant2$, the differential field $(k(X),d)$
is a finite algebraic extension of $(k(Y),\delta)$.

\smallskip

Our cyclotomic derivation $d$ is the  Jouanolou derivation with $s=1$,
and the cyclotomic derivation $\Delta$ is the factorisable derivation of $k[Y]$
associated with $d$.
In this case $s=1$, the differential field $(k(X),d)$ is no longer a finite algebraic extension of $(k(Y),\delta)$;
the relations between $d$ and $\Delta$ are thus more complicated.

We present some algebraic descriptions of the domains
$k[X]^d$, $k[Y]^\Delta$, and the fields $k(X)^d$, $k(Y)^\Delta$.
Note that these rings are nontrivial. The cyclic determinant
$$
w=\left|
\begin{array}{cccc}
x_0&x_1&\cdots&x_{n-1}\\
x_{n-1}&x_0&\cdots&x_{n-2}\\
\vdots&\vdots&&\vdots\\
x_2&x_3&\cdots&x_0\\
\end{array}
\right|
$$
is a polynomial belonging to $k[X]^d$, and
the product $y_0y_1\cdots y_{n-1}$ belongs to $k[Y]^\Delta$.
In this paper we prove, among others,
that $k(X)^d$ is a field of rational functions over $k$
in $n-\f(n)$ variables, where $\f$ is the Euler totient function (Theorem~\ref{DdTw}),
and that
$k[X]^d$ is a polynomial ring over $k$
if and only if $n$ is a power of a prime (Theorem~\ref{TwPol}).
The field $k(X)^d$ is in fact
the field of quotients of  $k[X]^d$ (Proposition~\ref{DdfPropB}).
We  denote by $\xi(n)$ the  sum $\sum_{p\mid n}\frac{n}{p}$,
where $p$ runs through all prime divisors of $n$, and we prove
that the number of the minimal set of generators of $k[X]^d$
is equal to $\xi(n)$
if and only if $n$ has at most two prime divisors (Corollary~\ref{MinWnFin}).
In particular,
if  $n=p^i q^j$, where $p\neq q$ are primes and $i,j$ are positive integers,
then the minimal number of generators of $k[X]^d$
is equal to $\xi(n)=p^{i-1}q^{j-1}(p+q)$ (Corollary~\ref{DwaWn}).

The ring of constants $k[Y]^\Delta$ is always equal to $k[v]$,
where $v=y_0y_1\dots,y_{n-1}$ (Theorem~\ref{TwDelta})
and, if $n$ is prime, then  $k(Y)^\Delta=k(v)$ (Theorem~\ref{ATtwPrime}).
If $n=p^s$, where $p$ is a prime and $s\geqslant 2$, then
$
k(Y)^\Delta=k(v, f_1,\dots,f_{m-1})
$
with $m=p^{s-1}$, where  $f_1,\dots,f_{m-1}\in k(Y)$
are homogeneous rational functions such that
$v,f_1,\dots,f_{m-1}$ are algebraically independent over $k$ (Theorem~\ref{RPStw}).
A similar theorem we prove for $n=p q$ (Theorem~\ref{PQtw}).

In our proofs we use classical properties of cyclotomic polynomials, and
an important role play some results (\cite{Lam}, \cite{Lam1}, \cite{Ste}, \cite{Ste1} and others)
on vanishing sums of roots of unity.

\section{Notations  and preparatory facts}\label{SecNP}

\initpar
We denote by $\mathbb Z_n$ the ring $\mathbb Z/n\mathbb Z$, and by $\mathbb Z_n^\ast$ the multiplicative
group of $\mathbb Z_n$.  The indexes of the variables $x_0,\dots,x_{n-1}$
and $y_0,\dots,y_{n-1}$ are elements
of $\mathbb Z_n$. This means, in particular,  that if $i,j$ are integers, then
$x_i=x_j\iff \MOD{i}{j}{n}$.
Throughout this paper $\E$ is a primitive $n$-th root of unity,
and  we assume that $\E\in k$. The letters $\varrho$ and $\tau$ we book
for two $k$-automorphisms of the field $k(X)=k\left(x_0,\dots,x_{n-1}\right)$,
defined by
$$
\varrho(x_j)=x_{j+1},\quad \tau(x_j)=\E^jx_j
\quad
\mbox{for all}\quad j\in\mathbb Z_n.
$$

We denote by $u_0,u_1,\dots,u_{n-1}$ the  linear forms in
$k[X]=k[x_0,\dots,x_{n-1}]$,
defined by
$$
u_j=\sum_{i=0}^{n-1} \left(\E^j\right)^ix_i,\quad\mbox{for}\quad j\in\mathbb Z_n.
$$
If $r$ is an integer and $n\nmid r$, then the sum
$\sum_{j=0}^{n-1}\left(\E^r\right)^j$ is equal to $0$, and in the other case,
when $n\mid r$, this sum is equal to $n$.
As a consequence of this fact we obtain, that
$$
x_i=\frac{1}{n}\sum_{j=0}^{n-1}\left(\E^{-i}\right)^j u_j\quad\mbox{for all}\quad i\in\mathbb Z_n.
$$
Thus, $k[X]=k[u_0,\dots,u_{n-1}]$, \  $k(X)=k(u_0,\dots,u_{n-1})$,
and the forms $u_0,\dots,u_{n-1}$ are algebraically independent over $k$.
Moreover, it is easy to check the following equalities.

\blem{NVlemC}
\ $\tau(u_j)=u_{j+1},\quad \varrho(u_j)=\E^{-j}u_j$
 \ for all $j\in\mathbb Z_n$.
\elem

For every sequence
$\alpha=(\alpha_0,\alpha_1,\dots,\alpha_{n-1})$, of integers,
we  denote by $H_\alpha(t)$  the polynomial in $\mathbb Z[t]$ defined by
$$
H_\alpha(t)=\alpha_0+\alpha_1t^1+\alpha_2t^2+\dots+\alpha_{n-1}t^{n-1}.
$$
An important role in our paper  play two subsets of $\mathbb Z^n$
which we denote by ${\cal G}_n $ and ${\cal M}_n $.
The first subset ${\cal G}_n$ is the set of all sequences
$\alpha=(\alpha_0,\dots,\alpha_{n-1})\in\mathbb Z^n$ such that
$
\alpha_0+\alpha_1\E^1+\alpha_2\E^2+\dots+\alpha_{n-1}\E^{n-1}=0.
$
The second  subset ${\cal M}_n $ is the set of all such sequences
$\alpha=(\alpha_0,\dots,\alpha_{n-1})$
which belong to ${\cal G}_n $ and the integers  $\alpha_0,\dots,\alpha_{n-1}$
are nonnegative, that is, they belong to the set of natural numbers
$\mathbb N=\{0,1,2,\dots\}$.
Let us remember:
$$
{\cal G}_n =\left\{\alpha\in\mathbb Z^n; \ H_\alpha(\E)=0\right\}, \quad
{\cal M}_n =\left\{\alpha\in\mathbb N^n; \ H_\alpha(\E)=0\right\}={\cal G}_n \cap\mathbb N^n.
$$
If $\alpha,\beta\in {\cal G}_n $, then of course $\alpha\pm\beta\in {\cal G}_n $,
and if $\alpha,\beta\in {\cal M}_n $, then $\alpha+\beta\in {\cal M}_n $.
Thus ${\cal G}_n $ is an abelian group, and ${\cal M}_n $ is an abelian monoid
with zero $0=(0,\dots,0)$.

\medskip

The primitive $n$-th root  $\E$ is an algebraic element over $\mathbb Q$, and its
monic minimal polynomial is equal to
the $n$-th cyclotomic polynomial $\Phi_n(t)$.
Recall  (see for example: \cite{Nagl},  \cite{Lang})
that $\Phi_n(t)$ is a monic irreducible polynomial with integer coefficients
of degree $\f(n)$, where $\f$ is the Euler totient function. This implies that we have the
following proposition.

\bprop{NVpropG}
Let $\alpha\in\mathbb Z^n$. Then
$\alpha\in {\cal G}_n $ if and only if there exists a polynomial $F(t)\in\mathbb Z[t]$
such that
$
H_\alpha(t)=F(t)\Phi_n(t).
$
\eprop

Put $e_0=(1,0,0,\dots,0)$, $e_1=(0,1,0,\dots,0), \ \dots, \ e_{n-1}=(0,0,\dots,0,1)$,
and let
$
e=\sum_{i=0}^{n-1}e_i=(1,1,\dots,1).
$
Since $\sum_{i=0}^{n-1}\E^i=0$, the element $e$ belongs to ${\cal M}_n $.

\bprop{NVpropMM}
If $\alpha\in {\cal G}_n $, then there exist $\beta,\gamma\in {\cal M}_n $ such that
$\alpha=\beta-\gamma$.
\eprop

\begin{dowod}
Let $\alpha=(\alpha_0,\dots,\alpha_{n-1})\in {\cal G}_n $, and
let $r=\min\{\alpha_0,\dots,\alpha_{n-1}\}$.
If $r\geqslant0$, then $\alpha\in {\cal M}_n $ and then $\alpha=\beta-\gamma$,
where $\beta=\alpha$, $\gamma=0$.
Assume that $r=-s$, where $1\leqslant s\in\mathbb N$.
Put $\beta=\alpha+se$ and $\gamma=se$. Then $\beta,\gamma\in {\cal M}_n $, and $\alpha=\beta-\gamma$.
\end{dowod}

\medskip

The monoid ${\cal M}_n $ has an order $\geqslant$.
If  $\alpha,\beta\in {\cal G}_n $, the we write $\alpha\geqslant\beta$,
if $\alpha-\beta\in\mathbb N^n$,
that is, $\alpha\geqslant\beta\iff$ there exists $\gamma\in {\cal M}_n $ such that
$\alpha=\beta+\gamma$.
In particular, $\alpha\geqslant0$ for any $\alpha\in {\cal M}_n $.
It is clear that the relation $\geqslant$ is
reflexive, transitive and antisymmetric. Thus ${\cal M}_n $ is a poset with respect to $\geqslant$.

\bprop{NVpropPos}
The poset ${\cal M}_n $ is artinian, that is,
if
$
\alpha^{(1)}\geqslant \alpha^{(2)}\geqslant \alpha^{(3)}\geqslant \dots
$
is a sequence of elements from ${\cal M}_n $, then there exists an integer $s$
such that $\alpha^{(j)}=\alpha^{(j+1)}$ for all $j\geqslant s$.
\eprop

\begin{dowod}
Given an element $\alpha=(\alpha_0,\dots,\alpha_{n-1})\in {\cal M}_n $,
we put $|\alpha|=\alpha_0+\dots+\alpha_{n-1}$.
Observe that if $\alpha,\beta\in {\cal M}_n $ and $\alpha>\beta$, then $|\alpha|>|\beta|$.
Suppose that there exists an infinite sequence
$
\alpha^{(1)}> \alpha^{(2)}> \alpha^{(3)}> \dots
$
of elements from ${\cal M}_n $, and let $s=\left|\alpha^{(1)}\right|$.
Then we have an infinite sequence $s>|\alpha^{(2)}|>|\alpha^{(2)}|>\dots
\geqslant 0$,
of natural numbers; a contradiction.
\end{dowod}

Let $\alpha\in {\cal M}_n $. We say that $\alpha$ is a {\it minimal element} of ${\cal M}_n $, if
$\alpha\neq0$ and there is no $\beta\in {\cal M}_n $ such that $\beta\neq0$ and $\beta<\alpha$.
Equivalently, $\alpha$ is a minimal element of ${\cal M}_n $, if $\alpha\neq0$ and $\alpha$ is not a sum
of two nonzero elements of ${\cal M}_n $.
It follows from Proposition \ref{NVpropPos} that for any $0\neq \alpha\in {\cal M}_n $
there exists a minimal element $\beta$ such that $\beta\leqslant\alpha$.
Moreover, every nonzero element of ${\cal M}_n $ is a finite sum of minimal elements.

\bprop{NVpropMin}
The set of all minimal elements of ${\cal M}_n $ is finite.
\eprop

\begin{dowod}
To deduce this result from Proposition \ref{NVpropPos}, Dikson's Lemma could be used : in any subset $\cal N$
of $\mathbb N^n$ there exists a finite number of elements $\{e^{(1)},\cdots,e^{(s)}\}$
such that ${\cal N} \subseteq \bigcup \left(e^{(j)} + \mathbb N^n\right)$.

It is simpler to use classical noetherian arguments.
Consider the polynomial ring $R=\mathbb Z[z_0,\dots,z_{n-1}]$.
If $\alpha=(\alpha_0,\dots,\alpha_{n-1})$ is an element from ${\cal M}_n $, then we denote
by $z^\alpha$ the monomial
$z_0^{\alpha_0}z_1^{\alpha_1}\cdots z_{n_1}^{\alpha_{n-1}}$.
Let ${\cal S}$ be the set of all minimal elements of ${\cal M}_n $, and consider
the ideal $A$ of $R$ generated by all elements of the form $z^\alpha$
with $\alpha\in{\cal S}$.
Since $R$ is noetherian, $A$ is finitely generated; there exist $\alpha^{(1)},\dots,\alpha^{(r)} \in {\cal S}$
such that
$
A=\left(z^{\alpha^{(1)}}, \ \dots, \ z^{\alpha^{(r)}}\right).
$
Let $\alpha$ be an arbitrary element from ${\cal S}$. Then $z^\alpha\in A$,
and then
there exist $j\in\{1,\dots,r\}$ and  $\gamma\in \mathbb N^n$ such
that $z^\alpha=z^\gamma\cdot z^{\alpha^{(j)}}=z^{\gamma+\alpha^{(j)}}$.
This implies that
$\alpha=\gamma+\alpha^{(j)}.$
Observe that $\gamma=\alpha-\alpha^{(j)}\in {\cal G}_n \cap\mathbb N^n$, and
${\cal G}_n \cap\mathbb N^n={\cal M}_n $, so $\gamma$ belongs to  ${\cal M}_n $.
But $\alpha$ is minimal, so $\gamma=0$, and consequently $\alpha=\alpha^{(j)}$.
This means that ${\cal S}$ is a finite set equal to
$\left\{\alpha^{(1)},\dots,\alpha^{(r)}\right\}$.
\end{dowod}

We denote by $\zeta$, the rotation of $\mathbb Z^n$ given by
$$
\zeta(\alpha)=\left(\alpha_{n-1},\alpha_0,\alpha_1,\dots,\alpha_{n-2}\right),
$$
for $\alpha=\left(\alpha_0,\alpha_1,\dots,\alpha_{n-1}\right)\in\mathbb Z^n$.
We have for example: $\zeta(e_j)=e_{j+1}$ for all $j\in\mathbb Z_n$,
and $\zeta(e)=e$.
The mapping
$\zeta:\mathbb Z^n\to\mathbb Z^n$ is obviously an endomorphism of the $\mathbb Z$-module
$\mathbb Z^n$, and is one-to-one and onto.

\blem{NVpropRot}
Let $\alpha\in\mathbb Z^n$.
If $\alpha\in {\cal G}_n$, then $\zeta(\alpha)\in {\cal G}_n$.
If $\alpha\in {\cal M}_n$, then $\zeta(\alpha)\in {\cal M}_n$.
Moreover, $\alpha$ is a minimal element of ${\cal M}_n$
if and only if $\zeta(\alpha)$ is a minimal element of~${\cal M}_n$.
\elem

\begin{dowod}
Assume that $\alpha=(\alpha_0,\dots,\alpha_{n-1})\in{\cal G}_n$.
Then $\alpha_0+\alpha_1\E+\dots+\alpha_{n-1}\E^{n-1}=0$. Multiplying it by $\E$, we have
$0=\alpha_0\E+\alpha_1\E^2+\dots+\alpha_{n-1}\E^{n}$.
But $\E^n=1$, so
$
\alpha_{n-1}+\alpha_0\E+\alpha_1\E^2+\dots+\alpha_{n-2}\E^{n-2}=0,
$
and so $\zeta(\alpha)\in{\cal G}_n$.
This implies also, that if $\alpha\in{\cal M}_n$, then $\zeta(\alpha)\in{\cal M}_n$.

Assume now that $\alpha$ is a minimal element of ${\cal M}_n$
and suppose that $\zeta(\alpha)=\beta+\gamma$, for some
$\beta,\gamma\in{\cal M}_n$. Then we have
$
\alpha=\zeta^n(\alpha)=\zeta^{n-1}\left(\zeta(\alpha)\right)
=\zeta^{n-1}(\beta)+\zeta^{n-1}(\gamma)=\beta'+\gamma',
$
where $\beta'=\zeta^{n-1}(\beta)$ and $\gamma'=\zeta^{n-1}(\gamma)$
belong to ${\cal M}_n$.
Since $\alpha$ is minimal, $\beta'=0$ or $\gamma'=0$, and then $\beta=0$ or $\gamma=0$.
Thus if $\alpha$ is a minimal element of ${\cal M}_n$, then
$\zeta(\alpha)$ is also a minimal element of ${\cal M}_n$.
Moreover, if $\zeta(\alpha)$ is minimal, then $\alpha$ is minimal,
because $\alpha=\zeta^{n-1}\left(\zeta(\alpha)\right)$.
\end{dowod}

\section{The derivation d and its constants}\label{SecCd}

Let us recall that $d:k[X]\to k[X]$ is a derivation such that $d(x_j)=x_{j+1}$,
for $j\in\mathbb Z_n$.

\bprop{DdPropA}
For each  $j\in \mathbb Z_n$, the equality
$
d(u_j)=\E^{-j}u_j
$
holds.
\eprop

\medskip

{\bf Proof.}
$
\begin{array}[t]{lcl}
d(u_j)&=&d\left(\sum\limits_{i=0}^{n-1}\left(\E^j\right)^i x_i\right)
=\sum\limits_{i=0}^{n-1}\left(\E^j\right)^{i} x_{i+1}
=\sum\limits_{i=1}^n\left(\E^j\right)^{i-1}x_i\medskip\\
&=&\E^{-j}\sum\limits_{i=1}^n\left(\E^j\right)^{i}x_i
=\E^{-j}\sum\limits_{i=0}^{n-1}\left(\E^j\right)^{i}x_i=
\E^{-j}u_j.  \ \kon
\end{array}
$
\medskip\\
This means that $d$ is a diagonal derivation of the polynomial ring
$k[U]=k[u_0,\dots,u_{n-1}]$ which is equal to the ring $k[X]$.
It is known (see for example \cite{Now}) that
the algebra of constants of
every diagonal derivation of $k[U]=k[X]$ is finitely generated over $k$.
Therefore, $k[X]^d$ is finitely generated over $k$.
We would like to describe a minimal set of generators of the ring $k[X]^d$,
and a minimal set of generators of the field $k(X)^d$.

\medskip

If $\alpha=(\alpha_0,\dots,\alpha_{n-1})\in\mathbb Z^n$, then we denote by $u^\alpha$
the rational monomial $u_0^{\alpha_0}\cdots u_{n-1}^{\alpha_{n-1}}$.
Recall (see the previous section) that
$H_\alpha(t)$ is the polynomial $a_0+a_1t^1+\dots+a_{n-1}t^{n-1}$
belonging to $\mathbb Z[t]$.
As a consequence of Proposition \ref{DdPropA} we obtain

\bprop{DdPropB}
\ $d(u^\alpha)=H_\alpha(\E^{-1})u^\alpha$
for all $\alpha\in\mathbb Z^n$.
\eprop

Note that $\E^{-1}$ is also a primitive $n$-th root of unity.
Hence, by Proposition \ref{NVpropG}, we have the equivalence
$H_\alpha(\E^{-1})=0\iff H_\alpha(\E)=0$, and so, by the previous proposition,
we see that if
$\alpha\in\mathbb Z^n$, then $d(u^\alpha)=0\iff \alpha\in{\cal G}_n$,
and if $\alpha\in\mathbb N^n$, then $d(u^\alpha)=0\iff \alpha\in{\cal M}_n$.
Moreover, if
$
F=b_1 u^{\alpha^{(1)}}+\dots+b_r u^{\alpha^{(r)}},
$
where  $b_1,\dots,b_r\in k$ and
$\alpha^{(1)},\dots,\alpha^{(r)}$ are pairwise distinct elements of
$\mathbb N^n$,
then
$d(F)=0$ if and only if $d\left(b_i u^{\alpha^{(i)}}\right)=0$
for every  $i=1,\dots,r$.
Hence, $k[X]^d$ is generated over $k$ by all elements
of the form $u^\alpha$ with $\alpha\in{\cal M}_n$.
We know (see the previous section), that every nonzero element of ${\cal M}_n$
is a finite sum of minimal elements of ${\cal M}_n$.
Thus we have the following next proposition.

\bprop{DdPropD}
The ring of constants $k[X]^d$ is generated over $k$ by all the elements
of the form $u^\beta$, where $\beta$ is a minimal element of the monoid
${\cal M}_n$.
\eprop

In the next section we will prove some additional facts on the minimal number
of generators of the ring $k[X]^d$.
Now, let us look at the field $k(X)^d$.

\bprop{DdfPropA}
The field of constants $k(X)^d$ is generated over $k$ by all elements
of the form $u^\gamma$ with $\gamma\in{\cal G}_n$.
\eprop

\begin{dowod}
Let $L$ be the subfield of $k(X)$ generated over $k$ by all elements
of the form $u^\gamma$ with $\gamma\in{\cal G}_n$.
It is clear that $L\subseteq k(X)^d$. We will prove the
reverse inclusion.
Assume that $0\neq f\in k(X)^d$. Since $k(X)=k(U)$, we have $f=A/B$,
where $A,B$ are coprime polynomials in $k[U]$.
Put
$$
A=\sum_{\alpha\in S_1} a_\alpha u^\alpha,\quad B=\sum_{\beta\in S_2} b_\beta
u^\beta,
$$
where all $a_\alpha$, $b_\beta$ are nonzero elements of
$k$, and $S_1$, $S_2$ are
some subsets of $\mathbb N^n$.
Since $d(f)=0$, we have the equality
$Ad(B)=d(A)B$. But $A,B$ are relatively prime, so $d(A)=\lambda A$, $d(B)=\Lambda B$
for some $\lambda\in k[U]$.
Comparing degrees, we see that $\lambda\in k$.
Moreover, by Proposition \ref{DdPropB}, we deduce that
$d(u^\alpha)=\lambda u^\alpha$ for all $\alpha\in S_1$, and also
$d(u^\beta)=\lambda u^\beta$ for all $\beta\in S_2$.
This implies that if $\delta_1,\delta_2\in S_1\cup S_2$, then $d\left(u^{\delta_1-\delta_2}\right)=0$.
In fact,
$
d\left(u^{\delta_1-\delta_2}\right)
=d\left(\frac{u^{\delta_1}}{u^{\delta_2}}\right)
=\frac{1}{u^{2\delta_2}}\left(d(u^{\delta_1})u^{\delta_2}-
u^{\delta_1}d(u^{\delta_2}\right)
=\frac{1}{u^{2\delta_2}}\left(\lambda u^{\delta_1}u^{\delta_2}-\lambda u^{\delta_1}u^{\delta_2}\right)=0.
$
This means, that if $\delta_1,\delta_2\in S_1\cup S_2$, then $\delta_1-\delta_2\in{\cal G}_n$
Fix an element $\delta$ from $S_1\cup S_2$.
Then all $\alpha-\delta$, $\beta-\delta$ belong to ${\cal G}_n$, and we have
$$
f=\frac{A}{B}=\frac{\sum a_\alpha u^\alpha}{\sum b_\beta u^\beta}
=\frac{u^{-\delta}\sum a_\alpha u^\alpha}{u^{-\delta}\sum b_\beta u^\beta}
=\frac{\sum a_\alpha u^{\alpha-\delta}}{\sum b_\beta u^{\beta-\delta}},
$$
and hence, $f\in L$.
\end{dowod}

Let us recall (see Proposition \ref{NVpropMM})
that every element of the group ${\cal G}_n$
is a difference of two elements from the monoid ${\cal M}_n$.
Using this fact and the previous propositions we obtain

\bprop{DdfPropB}
The field $k(X)^d$ is the field of quotients
of the ring  $k[X]^d$.
\eprop

Now we will prove that $k(X)^d$ is a field of rational functions
over $k$, and its transcendental degree over $k$ is equal to $n-\f(n)$, where
$\f$ is the Euler totient function.
For this aim look at the cyclotomic polynomial $\Phi_n(t)$.
Assume that
$$
\Phi_n(t)=c_0+c_1t+\dots+c_{\f(n)} t^{\f(n)}.
$$
All the coefficients $c_0,\dots,c_{\f(n)}$ are integers, and  $a_0=a_{\f(n)}=1$.
Put $m=n-\f(n)$ and
$$
\gamma_0=\Big(c_0,c_1,\dots,c_{\f(n)}, \ \underbrace{0,\dots,0}_{m-1}\Big).
$$
Note that $\gamma_0\in\mathbb Z^n$, and $H_{\gamma_0}(t)=\Phi_n(t)$.
Consider the elements $\gamma_0,\gamma_1,\dots,\gamma_{m-1}$ defined by
$$
\gamma_j=\zeta^j(\gamma_0),\quad\mbox{for}\quad j=0,1,\dots,m-1.
$$
Observe that $H_{\gamma_j}(t)=\Phi_n(t)\cdot t^j$
for all $j\in\{0,\dots,m-1\}$.
Since $\Phi_n(\E)=0$, we have $H_{\gamma_j}(\E)=0$, and so,
the elements $\gamma_0,\dots,\gamma_{m-1}$ belong to  ${\cal G}_n$.

\blem{DdLemG}
The elements $\gamma_0,\dots,\gamma_{m-1}$ generate the group ${\cal G}_n$.
\elem

\begin{dowod}
Let $\alpha\in {\cal G}_n$.
It follows from Proposition \ref{NVpropG},
that $H_\alpha(t)=F(t)\Phi_n(t)$, for some $F(t)\in\mathbb Z[t]$.
Then obviously $\deg F(t)<m$.
Put $F(t)=b_0+b_t^1+\dots+b_{m-1}t^{m-1}$, with $b_0,\dots,b_{m-1}\in\mathbb Z$.
Then we have
$$
\begin{array}{lcl}
H_\alpha(t)&=&b_0\left(\Phi_n(t)t^0\right)+b_1\left(\Phi_n(t)t^1\right)+\dots
+b_{m-1}\left(\Phi_n(t)t^{m-1}\right)\smallskip\\
&=&b_0 H_{\gamma_0}(t)+\dots+b_{m-1} H_{\gamma_{m-1}}(t),
\end{array}
$$
and this implies that
$\alpha=b_0\gamma_0+b_1\gamma_1+\dots+b_{m-1}\gamma_{m-1}$.
\end{dowod}

Consider now the rational monomials $w_0,\dots,w_{m-1}$ defined by
$$
w_j=u^{\gamma_j}=u_{0+j}^{c_0}u_{1+j}^{c_1}u_{2+j}^{c_2}\cdots u_{\f(n)+j}^{c_{\f(n)}}
$$
for  $j=0,1,\dots,m-1$,
where $m=n-\f(n)$.
Each $w_j$ is a rational monomial with respect
to  $u_0,\dots,u_{n-1}$
of the same degree equals to $\Phi_n(1)=c_0+c_1+\dots+c_{\f(n)}$.
It is known (see for example \cite{Lang}) that $\Phi_n(1)=p$ if $n$
is power of a prime number~$p$, and $\Phi_n(1)=1$ in all other cases.
As each $u_j$ is a homogeneous polynomial in $k[X]$ of degree $1$,
we have:

\bprop{DdHom}
The  elements $w_0,\dots,w_{m-1}$ are homogeneous rational functions
with respect to variables $x_0,\dots,x_{n-1}$, of the same degree $r$.
If $n$ is a power of a prime number~$p$, then $r=p$, and $r=1$ in all other cases.
\eprop

As an immediate consequence of Lemma \ref{DdLemG} and Proposition \ref{DdfPropA},
we obtain the equality $k(X)^d=k(w_0,\dots,w_{n-1})$.

\blem{DdLemAN}
The elements $w_0,\dots,w_{m-1}$ are algebraically independent over $k$.
\elem

\begin{dowod}
Let $A$ be the $n\times m$ Jacobi matrix $\left[a_{ij}\right]$, where
$a_{ij}=\frac{\partial w_j}{\partial u_i}$ for
$i=0,1,\dots,n-1$, \ $j=0,1,\dots,m-1$.
It is enough  to show that rank$(A)=m$ (see for example \cite{Jacob}).
Observe that
$\frac{\partial w_0}{\partial u_0}=c_0 u_0^{c_0-1}u_1^{c_1}\cdots u_{\f(n)}^{c_{\f(n)}}\neq0$
(because $c_0=1$), and $\frac{\partial w_j}{\partial u_0}=0$ for $j\geqslant1$.
Moreover, $\frac{\partial w_1}{\partial u_1}\neq0$ and
$\frac{\partial w_j}{\partial u_1}=0$ for $j\geqslant2$, and in general,
$\frac{\partial w_i}{\partial u_i}\neq 0$ and $\frac{\partial u_j}{\partial u_i}=0$
for all $i,j=0,\dots,m-1$ with $j>i$.
This means, that the upper $m\times m$ matrix of $A$ is a triangular matrix
with a nonzero determinant. Therefore, rank$(A)=m$.
\end{dowod}

Thus, we proved the following theorem.

\btw{DdTw}
The field of constants $k(X)^d$ is a field of rational functions over $k$
and its transcendental degree over $k$ is equal to
$m=n-\f(n)$, where $\f$ is the Euler totient function.
More precisely,
$$
k(X)^d=k\Big(w_0,\dots,w_{m-1}\Big),
$$
where the elements $w_0,\dots,w_{m-1}$ are as above.
\etw

Now we will describe all constants of $d$ which are
homogeneous rational functions of degree zero.
Let us recall that a nonzero polynomial $F$ is homogeneous of degree $r$,
if all its monomials are of the same degree $r$.
We assume that the zero polynomial
is homogeneous of arbitrary degree.
Homogeneous polynomials are also homogeneous rational functions,
which (in characteristic zero) are defined
in the following way.
Let $f=f(x_0,\dots,x_{n-1})\in k(X)$
We say that $f$ is {\it homogeneous} of degree $s\in\mathbb Z$,
if in the field $k(t,x_0,\dots,x_{n-1})$ the equality
$
f(tx_0,tx_1,\dots,tx_{n-1})=t^s\cdot f(x_0,\dots,x_{n-1})$ holds
It is easy to prove (see for example \cite{Now} Proposition 2.1.3)
the following equivalent formulations of homogeneous rational functions.

\bprop{HomHom}
Let $F,G$ be nonzero coprime polynomials in $k[X]$ and let $f=F/G$.
Let $s\in\mathbb Z$.
The following conditions are equivalent.


$(1)$ The rational function $f$ is homogeneous of degree $s$.


$(2)$ The polynomials $F$, $G$ are homogeneous of degrees $p$ and $q$,
respectively, where $s=p-q$.


$(3)$ \ $
x_0\frac{\partial f}{\partial x_0}+\dots+x_{n-1}\frac{\partial f}{\partial x_{n-1}}=sf$.
\eprop
Equality $(3)$ is called the {\it Euler formula}.
In this paper we denote by $E$ the {\it Euler derivation} of $k(X)$,
that is, $E$ is a derivation of $k(X)$  defined by
$E(x_j)=x_j$ for all $j\in\mathbb Z_n$.
As usually, we denote by $k(X)^E$ the field of constants of $E$.
Observe that, by Proposition \ref{HomHom}, a rational function $f\in k(X)$
belongs to $k(X)^E$ if and only if $f$ is homogeneous of degree zero.
In particular, the set of all homogeneous rational functions of degree zero
is a subfield of $k(X)$.
It is obvious that the quotients $\frac{x_1}{x_0},\dots,\frac{x_{n-1}}{x_0}$
belong to $k(X)^E$, and they are algebraically independent over $k$.
Moreover, $k(X)^E=k(\frac{x_1}{x_0},\dots,\frac{x_{n-1}}{x_0})$.
Therefore, $k(X)^E$ is a field of rational functions over $k$, and its transcendence degree over $k$
is equal to $n-1$.
Put
$q_j=\frac{x_{j+1}}{x_j}$ for all $j\in\mathbb Z_n$.
In particular, $q_{n-1}=\frac{x_0}{x_{n-1}}$.
The elements $q_0,\dots,q_{n-1}$ belong to $k(X)^E$ and moreover,
$\frac{x_j}{x_0}=q_0q_1\cdots q_{j-1}$ for $j=1,\dots,n-1$.
Thus we have the following equality.

\bprop{HomPropA}
\ $
k(X)^E=k\left(\frac{x_1}{x_0},\frac{x_2}{x_1},\dots,\frac{x_{n-1}}{x_{n-2}},\frac{x_0}{x_{n-1}}\right)$.
\eprop

Now consider the field
$
k(X)^{d,E}=k(X)^d\cap k(X)^E.
$

\blem{HomLemDer}
Let $d_1,d_2:k(X)\to k(X)$ be two derivations. Assume that
$K(X)^{d_1}=k(c,b_1,\dots,b_s)$, where $c,b_1,\dots,b_s$
are algebraically independent over $k$ elements from $k(X)$
such that $d_2(b_1)=\dots=d_2(b_s)=0$ and $d_2(c)\neq0$.
Then
$
k(X)^{d_1}\cap k(X)^{d_2}=k(b_1,\dots,b_s).
$
\elem

\begin{dowod}
Put $L=k(b_1,\dots,b_s)$. Observe that $k(X)^{d_1}=L(c)$,
and $c$ is transcendental over $L$.
Let $0\neq f\in k(X)^{d_1}\cap k(X)^{d_2}$.
Then $f=\frac{F(c)}{G(c)}$, where $F(t), G(t)$ are coprime polynomials
in $L[t]$. We have: $d_2(F(c))=F'(c)d_2(c)$, \ $d_2(G(c))=G'(c)d_2(c)$,
where $F'(t),G'(t)$ are  derivatives of $F(t), G(t)$, respectively.
Since $d_2(f)=0$, we have
$$
0=d_2(F(c))G(c)-d_2(G(c))F(c)=\Big(F'(c)G(c)-G'(c)F(c)\Big)d_2(c),
$$
and so, $(F'G-G'F)(c)=0$, because $d_2(c)\neq0$.
Since $c$ is transcendental over $L$, we obtain the equality
$
F'(t)G(t)=G'(t)F(t)
$
in $L[t]$, which implies that $F(t)$ divides $F'(t)$ and $G(t)$ divides $G'(t)$
(because  $F(t), G(t)$ are relatively prime), and comparing degrees
we deduce that $F'(t)=G'(t)=0$, that is, $F(t)\in L$ and $G(t)\in L$.
Thus the elements $F(c)$, $G(c)$
belong to $L$ and so, $f=\frac{F(c)}{G(c)}$
belongs to $L$. Therefore, $k(X)^{d_1}\cap k(X)^{d_2}\subseteq L$.
The reverse inclusion is obvious.
\end{dowod}

Let us return to the rational functions $w_0,\dots,w_{m-1}$.
We know (see Proposition \ref{DdHom}) that they are homogeneous
of the same degree. Put:
$d_1=d$, \ $d_2=E$, \ $c=w_0$ and
$b_j=\frac{w_j}{w_0}$  for $j=1,\dots,m-1$,
Then, as a consequence of Lemma \ref{HomLemDer}.
we obtain the following proposition.

\bprop{HomDerZero}
\ $
k(X)^{d,E}=k\left(\frac{w_1}{w_0},\dots,\frac{w_{m-1}}{w_0}\right)$.
\eprop
Since  $w_0,\dots,w_{m-1}$ are algebraically independent over $k$
(see Lemma \ref{DdLemAN}), the quotients $\frac{w_1}{w_0},\dots,\frac{w_{m-1}}{w_0}$
are also algebraically independent over $k$.
Thus, $k(X)^{d,E}$ is a field of rational functions
and its transcendental degree over $k$ is equal to
$n-\f(n)-1$, where $\f$ is the Euler totient function.
In particular, if  $n$ is prime, then $n-\f(n)-1=0$ and we obtain:

\bwn{HomWnPrime}
\ $k(X)^{d,E}=k\iff n$ is a prime number.
\ewn

\section{Numbers of minimal elements}\label{SecMin}

\initpar
Let  ${\cal F}$ be the set of all the minimal elements of the
monoid ${\cal M}_n$, and denote by $\nu(n)$ the cardinality of ${\cal F}$.
We know, by
Proposition~\ref{NVpropMin}, that $\nu(n)<\infty$.
We also  know (see Proposition \ref{DdPropD}) that the ring
$k[X]^d$ is generated over $k$
by all the elements of the form $u^\beta$, where $\beta\in{\cal F}$.
But $k[X]$
is equal to the polynomial ring $k[U]=k[u_0,\dots,u_{n-1}]$,
so $k[X]^d$ is generated over $k$ by a finite set of monomials
with respect to the variables $u_0,\dots,u_{n-1}$.

It is clear that if $\beta,\gamma$ are distinct elements from ${\cal F}$, then
$u^\beta\nmid u^\gamma$ and $u^\gamma\nmid u^\beta$.
This implies that no monomial $u^{\beta}, \beta\in{\cal F}$
belongs to the algebra generated by other $u^{\gamma}, \gamma\in{\cal F}, u^\gamma\nmid u^\beta$.
Thus, $\left\{u^\beta; \ \beta\in{\cal F}\right\}$
is a minimal set of generators of $k[X]^d$.

Moreover, $\left\{u^\beta; \ \beta\in{\cal F}\right\}$ is a set on generators of $k[X]^d$
with the minimal number of elements according to the following proposition.

\bprop{BGprop}
Let $f_1,\dots,f_s$ be polynomials in $k[X]$.
If $k[X]^d=k[f_1,\dots,f_s]$, then $s\geqslant \nu(n)$.
\eprop

\begin{dowod}
As the $u^\beta$ are monomials in the $u$'s, they constitute a Gr\"obner base for the ideal $I$ generated
in $k[X]$ by $k[X]^d$. This basis is minimal for any admissible order, for example the lexicographical one.

Making a head reduction of the $f_i$, a new head-reduced system of generators appears,
maybe with less than $s$  elements. Thus, without loss of generality, we can suppose that
the system $(f_1,\dots,f_s)$ is head-reduced, which means that the leading monomial
of one $f_i$ does not belong to the multiplicative monoid generated by the other leading monomials.

The leading monomials of the various $f_i$ are $u^\alpha$ for some $\alpha \in {\cal M}_n$.

The exponents $\alpha$ are minimal in the sub-monoid they generate, but this sub-monoid
has to be ${\cal M}_n$ itself.
\end{dowod}

\bigskip

In this section we prove, among others,
 that $k[X]^d$ is a polynomial ring over $k$
if and only if $n$ is a power of a prime number.
Moreover, we present some additional properties of the number $\nu(n)$,
which are consequences of
known results
on vanishing sums of roots of unity;
see for example \cite{Lam1}, \cite{Sat}, \cite{Ste} and \cite{Ste1},
where many interesting facts and references on this subject can be found.

\smallskip

We  denote by $\xi(n)$ the  sum $\sum_{p\mid n}\frac{n}{p}$,
where $p$ runs through all prime divisors of $n$.
Note that if  $a,b$ are
positive coprime integers, then $\xi(ab)=a\xi(b)+\xi(a)b$.

\medskip

First we  show that the computation of $\nu(n)$ can be reduced to the case
when $n$ is square-free.
For this aim
let us denote by $n_0$ the largest square-free factor of $n$, and by~$n'$ the integer $n/n_0$.
Then $\f(n)=n'\f(n_0)$, \ $\Phi_n(t)=\Phi_{n_0}\left(t^{n'}\right)$ (see for example \cite{Nagl}),
and  $\displaystyle \xi(n)=n'\xi(n_0)$.

\smallskip

Assume now that $n=mc$, where $m\geqslant2$, $c\geqslant2$ are integers.
For a given  sequence $\gamma=(\gamma_0,\dots,\gamma_{m-1})\in\mathbb Z^m$,
consider the sequence
$$
\overline{\gamma}=\Big(\gamma_0,\underbrace{0,\dots,0}_{c-1},\gamma_1,
\underbrace{0,\dots,0}_{c-1},\dots,\gamma_{m-1},\underbrace{0,\dots,0}_{c-1}\Big).
$$
This sequence is an element of $\mathbb Z^n$, and it is easy
to prove the following lemma.

\blem{MinLemB}
$\overline{\gamma}\in {\cal G}_n\iff \gamma\in{\cal G}_m$,
and $\overline{\gamma}\in {\cal M}_n\iff \gamma\in{\cal M}_n$.
Moreover,
$\overline{\gamma}$ is a minimal element of ${\cal M}_n\iff \gamma$
is a minimal element of ${\cal M}_m$.
\elem
Using the above notations, we have:
\bprop{MinPropA}
\ $\nu(n)=n'\nu(n_0)$, for all $n\geqslant3$.
\eprop

\begin{dowod}
If $n'=1$ then this is clear. Assume that $n'\geqslant2$.
Let $\alpha=(\alpha_0,\dots,\alpha_{n-1})$
be an element of ${\cal M}_n$. For every $j\in\{0,1,\dots,n'-1\}$, let us
denote:
$$
f_j(t)=\sum_{i=0}^{n_0-1} \alpha_{in'+j}t^{in'+j}=t^j \sum_{i=0}^{n_0-1} \alpha_{in'+j}t^{in'},
\quad \beta_j=\Big(\alpha_{0n'+j}, \ \alpha_{1n'+j}, \ \dots, \ \alpha_{(n_0-1)n'+j}\Big),
$$
 Note that $f_j(t)\in\mathbb Z[t]$
and $\beta_j\in\mathbb N^{n_0}$.
Consider the elements $\overline{\beta_0}, \  \overline{\beta_1}, \
\dots, \ \overline{\beta_{n'-1}}$, introduced before Lemma \ref{MinLemB} for $m=n_0$ and $c=n'$.
Observe that
$$
\alpha=\overline{\beta_0}+\zeta(\overline{\beta_1})
+\zeta^2(\overline{\beta_2})+\dots+\zeta^{n'-1}(\overline{\beta_{n'-1}})
\leqno{(\ast)}
$$
where $\zeta$ is the rotation of $\mathbb Z^n$, as in Section \ref{SecNP}.
Denote also by $f(t)$ the polynomial
$H_\alpha(t)=\alpha_0+\alpha_1t+\dots+\alpha_{n-1}t^{n-1}$,
that is, $f(t)=\sum_{j=0}^{n'-1} f_j(t)$.
It follows from  Proposition~\ref{NVpropG}, that $f(t)=g(t)\Phi_n(t)$
for some  $g(t)\in \mathbb Z[t]$.

\smallskip

For every $j\in\{0,1,\dots,n'-1\}$, denote by
$A_j$ the set of  polynomials
$F(t)\in \mathbb Z[t]$ such that the degrees of all nonzero monomials of $F(t)$
are congruent to $j$ modulo $n'$.  We assume that the zero polynomial
also belongs to $A_j$. It is clear that each $A_j$ is a $\mathbb Z$-module,
$A_iA_j\subseteq A_{i+j}$ for $i,j\in\mathbb Z_{n'}$,
and $\mathbb Z[t]=\bigoplus_{j\in\mathbb Z_{n'}} A_j$.
Thus, we  have a gradation on $\mathbb Z[t]$ with respect to $\mathbb Z_{n'}$.
We will say that it is the {\it $n'$-gradation}, and the decompositions of
polynomials with respect to this gradations we will call the {\it $n'$-decompositions}.

\smallskip

Let $g(t)=g_0(t)+g_1(t)+\dots+g_{n'-1}(t)$ be the $n'$-decomposition of $g(t)$;
each $g_j(t)$ belongs to $A_j$. Since $\Phi_n(t)=\Phi_{n_0}(t^{n'})$,
$\Phi_n(t)\in A_0$ and
$$
f(t)=g_0(t)\Phi_n(t)+g_1(t)\Phi_n(t)+\dots+g_{n'-1}(t)\Phi_n(t),
$$
is the $n'$-decomposition of $f(t)$.
But the previous equality $f(t)=\sum f_j(t)$ is also the
$n'$-decomposition of $f(t)$, so we have
$f_j(t)=g_j(t)\Phi_n(t)$ for all $j\in\mathbb Z_{n'}$.

Put $\eta=\E^{n'}$. Then $\eta$ is a primitive $n_0$-th root of unity
and, for every $j\in\mathbb Z_{n'}$,
$$
\sum_{i=0}^{n_0-1}\alpha_{in'+j}\eta^i=\E^{-j}f_j(\E)=\E^{-j}g_j(\E)\Phi_n(\E)
=\E^{-j}g_j(\E)\cdot0=0.
$$
This means that each $ \beta_j$ is an element of ${\cal M}_{n_0}$.

\smallskip

Assume now that the above $\alpha$ is a minimal element of ${\cal M}_n$.
Then, by~$(\ast)$, we have $\alpha=\zeta^j(\overline{\beta_j})$
for some $j\in\{0,\dots,n'-1\}$. Then
$\overline{\beta_j}=\zeta^{n-j}(\alpha)$ and so,
$\overline{\beta_j}$ is (by Lemma \ref{NVpropRot})
a minimal element of ${\cal M}_n$, and this implies, by Lemma \ref{MinLemB},
that $\beta_j$ is a minimal element of ${\cal M}_{n_0}$.
Thus, every minimal element $\alpha$ of ${\cal M}_n$
is of the form $\alpha=\zeta^j(\overline{\beta})$,
where $j\in\{0,\dots,n'-1\}$
and $\beta$ is a minimal element of ${\cal M}_{n_0}$, and it is clear that
this presentation is unique.
This means, that $\nu(n)\leqslant n'\cdot \nu(n_0)$.

\smallskip

Assume now that  $\beta$ is a minimal element of ${\cal M}_{n_0}$.
Then we have
$n'$ pairwise distinct sequences
$
\overline{\beta}, \  \zeta(\overline{\beta}), \  \zeta^2(\overline{\beta}),
\ \dots, \ \zeta^{n'-1}(\overline{\beta}),
$
which are (by Lemmas \ref{NVpropRot} and \ref{MinLemB})
minimal elements of ${\cal M}_n$.
Hence, $\nu(n)\geqslant n'\cdot \nu(n_0)$. Therefore, $\nu(n)=n'\cdot\nu(n_0)$.
\end{dowod}

\medskip


If $p$ is prime, then $\nu(p)=1$; the constant sequence $e=(1,1,\dots,1)$
is a unique minimal element of ${\cal M}_p$. In this case
$k[X]^d$ is the polynomial ring
$k[w]$,
where $w=u_0\dots u_{p-1}$
is the cyclic determinant of the variables $x_0,\dots,x_{p-1}$
(see Introduction). In particular, if $p=3$, then
$k[x_0,x_1,x_2]^d=k[x_0^3+x_1^3+x_2^3-3x_0x_1x_2]$.
Using Proposition \ref{MinPropA} and its proof we obtain:

\bprop{MinPot}
Let  $n=p^{s}$, where $s\geqslant1$ and $p$ is a prime number.
Then $\nu(n)=\xi(n)=p^{s-1}$,
and the ring of constants  $k[X]^d$ is a polynomial ring over $k$ in $p^{s-1}$ variables.
\eprop

\medskip


Assume now that $p$ is a prime divisor of $n$.  Denote by  $n_p$ the integer
$n/p$, and consider the sequences
$$
E^{(p)}_i=\sum_{j=0}^{p-1} e_{i+jn_p},
$$
for $i=0,1,\dots,n_p-1$.  Recall that $e_0=(1,0,\dots,0), \ \dots, \ e_{n-1}=(0,0,\dots,0,1)$
are the basic elements of $\mathbb Z^n$. Observe that each $E^{(p)}_i$
is equal to $\zeta^i\left(E^{(p)}_0\right)$, where $\zeta$ is the rotation of~$\mathbb Z^n$.
Observe also that $E^{(p)}_0=\overline{e}$, where in this case $e=(1,1,\dots,1)\in \mathbb Z^p$
and $\overline{e}$ is the element of $\mathbb Z^n$ introduced before Lemma \ref{MinLemB}
for $m=p$ and $c=n_p$. But $e$ is a minimal element of ${\cal M}_p$, so we see, by Lemmas
\ref{MinLemB} and \ref{NVpropRot}, that each $E^{(p)}_i$ is a minimal element of ${\cal M}_n$.
We will say that such $E^{(p)}_i$ is a {\it standard} minimal element of ${\cal M}_n$.
It is clear that if $i,j\in\{0,1,\dots,n_p-1\}$ and $i\neq j$, then
$E^{(p)}_i\neq E^{(p)}_j$. Observe also that, for every $i$, we have
$\left|E^{(p)}_i\right|=p$. This implies, that if $p\neq q$ are prime divisors of $n$, then
$E^{(p)}_i\neq E^{(q)}_j$ for all $i\in\{0,\dots,n_p-1\}$, \ $j\in\{0,1,\dots,n_q-1\}$.
Assume that $p_1,\dots,p_s$ are all the prime divisors of $n$.
Then, by the above observations, the number of all standard minimal elements of ${\cal M}_n$
is equal to $n_{p_1}+\dots+n_{p_s}$, that is, it is equal to $\xi(n)$.
Hence, we proved the following proposition.

\bprop{MinPropNrw}
\ $\nu(n)\geqslant\xi(n)$,  for all $n\geqslant3$.
\eprop

For a proof of the next result we need the following lemma.

\blem{MinABC}
If $n$ is divisible by two distinct primes, then $\xi(n)+\f(n)>n$.
\elem

\begin{dowod}
Since $\xi(n)=n'\xi(n_0)$, $\f(n)=n'\f(n_0)$ and $n=n'n_0$
we may assume that $n$ is square-free.
Let $n=p_1\cdots p_s$, where $s\geqslant2$ and $p_1,\dots,p_s$ are distinct primes.
If $s=2$, then the equality is obvious.
Assume that $s\geqslant3$, and that the equality is true for $s-1$.
Put $p=p_s$, $m=p_1\cdots p_{s-1}$. Then $m$ is square-free, $n=mp$,
$\gcd(m,p)=1$, $\xi(m)+\f(m)>m$ and moreover, $\f(m)<m$.
Hence,
$\xi(n)+\f(n)=p\xi(m)+\xi(p)m+\f(p)\f(m)=p\xi(m)+m+(p-1)\f(m)>p\xi(m)+p\f(m)>pm=n$.
and hence, by an induction, $\xi(n)+\f(n)>n$.
\end{dowod}

\btw{TwPol}
The ring of constants $k[X]^d$ is a polynomial ring over $k$
if and only if $n$ is a power of a prime number.
\etw

\begin{dowod}
Assume that $n$ is
divisible by two distinct primes, and suppose that $k[X]^d$ is a polynomial
ring of the form $k[f_1,\dots,f_s]$, where $f_1,\dots,f_s\in k[X]$
are algebraically independent over $k$.
Then, by Proposition \ref{BGprop}, we have $s\geqslant \nu(n)$.
The polynomials  $f_1,\dots,f_s$ belong to the field $k(X)^d$,
and we know, by Theorem \ref{DdTw}, that the transcendental degree
of this field over $k$ is equal to $n-\f(n)$.
Hence, $s\leqslant n-\f(n)$. But $\nu(n)\geqslant \xi(n)$
(Proposition~\ref{MinPropNrw}) and $\xi(n)>n-\f(n)$ (Lemma~\ref{MinABC}),
so we have a contradiction: $s\geqslant \nu(n)\geqslant \xi(n)>n-f(n)$.
This means, that if $n$ is divisible by two distinct primes, then $k[X]^d$
is not a polynomial ring over $k$. Now this theorem follows
from Proposition~\ref{MinPot}.
\end{dowod}

It is well known (see for example \cite{BeS}) that all coefficients of
the cyclotomic polynomial $\Phi_n(t)$ are nonnegative if and only if
$n$ is a power of a prime. Thus, we proved that $k[X]^d$ is a polynomial ring over $k$
if and only if all coefficients of $\Phi_n(t)$ are nonnegative.

\medskip

In our next considerations we will apply the following
theorem of R\'edei, de Bruijn and Schoenberg.

\bbtw{\cite{Red}, \cite{deB}, \cite{Sch}}{TwRBS}
The standard minimal elements of ${\cal M}_n$ generate the group ${\cal G}_n$.
\etw

Known proofs of the above theorem used usually
techniques of group rings.
Lam and Leung \cite{Lam1} gave a new proof using induction and group-theoretic
techniques.


\medskip

Now, let us  assume that $n=pq$, where $p\neq q$ are primes.
In this case, Lam and Leung \cite{Lam1} proved that
$\nu(n)=p+q$.
We will give a new elementary proof of this fact.
Note that in this case $n_p=q$ and $n_q=p$,
Put $P_i=E^{(q)}_i$ for $i=0,1,\dots,p-1$, and $Q_j=E^{(p)}_j$ for $j=0,\dots,q-1$.
We have $p+q$
elements $P_0,\dots, P_{p-1}, \  Q_0, \dots, Q_{q-1}$,
which are the standard minimal elements of ${\cal M}_{pq}$.

\blem{MinPQlemC}
For every $\beta\in{\cal M}_{pq}$ there exist nonnegative integers
$a_0,\dots, a_{p-1}$, $b_0,\dots$, $b_{q-1}$ such that
$
\beta=a_0P_0+\dots+a_{p-1}P_{p-1}+b_0Q_0+\dots+b_{q-1}Q_{q-1}.
$
\elem

\begin{dowod}
Let $\beta\in{\cal M}_{pq}$. Then $\beta\in{\cal G}_{pq}$ and, by Theorem \ref{TwRBS},
we have an equality $\beta=\sum a_iP_i+\sum b_jQ_j$, for some integers
$a_0,\dots,a_{p-1},b_0,\dots,b_{q-1}$.
Since
$\sum_{i=0}^{p-1} P_i=e=\sum_{j=0}^{q-1} Q_j$,
we may assume that $b_{q-1}=0$.
Let us recall that  $P_i=\sum_{j=0}^{q-1}e_{jp+i}$ for $i=0,\dots,p-1$,
and
$Q_j=\sum_{i=0}^{p-1}e_{iq+j}$ for $j=0,\dots,q-1$. Thus, we have
$$
\beta=\sum_{i=0}^{p-1}\sum_{j=0}^{q-1}\Big(a_i e_{jp+i}+b_je_{iq+j}\Big).
\leqno{(1)}
$$
Every number $m$ from $\{0,1,\dots,pq-1\}$
has a unique presentation in the form $m=sp+r$
with $s\in\{0,\dots,q-1\}$, $r\in\{0,\dots,p-1\}$, and it has also a unique
presentation $m=s_1q+r_1$ with
$s_1\in\{0,\dots,p-1\}$, $r_1\in\{0,\dots,q-1\}$.
Hence, it follows from $(1)$ that
$$a_i+b_j\geqslant0\quad\mbox{for all}\quad
i\in\{0,\dots,p-1\}, \ j\in\{0,\dots,q-1\}.
\leqno{(2)}
$$
But $b_{q-1}=0$, so $a_i\geqslant0$ for all $i=0,\dots,p-1$.
If all the numbers $b_0,\dots,b_{q-2}$ are also nonnegative, then we are done.

Assume that among $b_0,\dots,b_{q-2}$ there exists a negative integer,
and consider the number
$b_s=\min\{b_0,\dots,b_{q-2}\}$. Then $s\in\{0,\dots,q-2\}$
and $-b_s>0$. Put $A=\{0,\dots,q-1\} \smallsetminus \{s\}$.
Using again the equality $\sum_{i=0}^{p-1} P_i=\sum_{j=0}^{q-1} Q_j$,
we have: $b_sQ_s=\sum_{i=0}^{p-1} b_sP_i+\sum_{j\in A} (-b_s)Q_j$.
Hence,
$$
\beta=\sum_{i=0}^{p-1}(a_i+b_s)P_i+\sum_{j\in A} (b_j-b_s)Q_j+(-b_s)Q_{q-1}.
$$
By $(2)$, each $a_i+b_s$ is nonnegative. Moreover $b_s\leqslant b_s$ for all $j\in A$, and $-b_s>0$.
Therefore, in the above presentation all the
coefficients are nonnegative integers.
\end{dowod}

\bbtw{\cite{Lam1}}{MinPQtw}
Let $n=p^iq^j$, where $p\neq q$ are primes and $i,j$ are positive integers.
Then
$
\nu(n)=\xi(n)=p^{i-1}q^{j-1}(p+q).
$
In other words, the monoid ${\cal M}_n$ has exactly
$p^{i-1}q^{j-1}(p+q)$ minimal elements, and all its minimal elements are standard.
\etw

\begin{dowod}
Let $n=pq$, and ${\cal B}=\{P_0,\dots,P_{p-1},Q_0,\dots,Q_{q-1}\}$.
We know
that every element of ${\cal B}$ is a standard minimal element
of ${\cal M}_{pq}$, and that
all these elements are pairwise distinct.
Moreover, it follows from Lemma \ref{MinPQlemC} that every $\beta\in{\cal M}_{pq}$,
which is a minimal element of ${\cal M}_{pq}$, belongs to  ${\cal B}$.
Hence,  $\nu(pq)=p+q=\xi(pq)$.
This implies, by the equality $\xi(n)=n'\xi(n_0)$ and Proposition~\ref{MinPropA},
that $\nu(n)=\xi(n)$ for all $n$ of the form $p^iq^j$.
\end{dowod}

As a consequence of Theorem \ref{MinPQtw} and Proposition \ref{BGprop}
we obtain:

\bwn{DwaWn}
Let $n=p^iq^j$, where $p\neq q$ are primes and $i,j$ are positive integers.
Then the minimal number of generators of the ring of constants $k[X]^d$
is equal to $\xi(n)=p^{i-1}q^{j-1}(p+q)$.
\ewn

\medskip

We already know
that if $n$ is divisible by at most two distinct primes,
then every minimal element of ${\cal M}_n$ is standard.
It is well known (see for example \cite{Lam1}, \cite{Ste1}, \cite{Sat})
that in all other cases always exist
nonstandard minimal elements.
For instance, Lam and Leung \cite{Lam1} proved
that if $n$ is divisible by three primes $p_1<p_2<p_3$, then
the equality $a_1a_2+a_3=0$, where $a_j=\sum_{i=1}^{p_1-1}\E^{in_{p_i}}$
for $j=1,2,3$, is of the form $H_\alpha(\E)=0$, where  $\alpha$
is a nonstandard minimal element of ${\cal M}_n$.
There are also other examples. Assume that $n=p_1\cdots p_s$, where
$p_1,\dots,p_s$ are distinct primes.
and denote by  $U$ the set of all numbers from $\{1,2,\dots,n-1\}$
which are relatively prime to $n$.
If $s\geqslant3$ is odd, then
$$
\gamma=e_0+\sum_{u\in U} e_u.
$$
is a nonstandard minimal element of ${\cal M}_n$.
This element $\gamma$ belongs to ${\cal M}_n$,
because
the sum of all primitive $n$-th roots of unity is equal to $\mu(n)$,
where $\mu$ is the M\"obius function (see for example \cite{L-N}, \cite{Mot}).
The minimality of $\gamma$ follows from the known fact (see for example \cite{Con})
that if $n$ is square-free, then all the primitive $n$-th roots of unity
form a basis of $\mathbb Q(\E)$ over $\mathbb Q$.
Observe also that $|\gamma|=\f(n)+1\neq p_i$ for all $i=1,\dots,s$,
so $\gamma$ is nonstandard.

If $s\geqslant 4$ is even, then
put $p=p_s$ ,  $n'=p_1\cdots p_{s-1}$, and let $U'$
the set of all numbers from  $\{1,2,\dots,n'-1\}$
which are relatively prime to $n'$.
Then  $\E^p$ is a primitive $n'$-th root of unity and, using similar arguments,
we see that
$$
\gamma'=e_0+\sum_{v\in U'} e_{vp}.
$$
is a nonstandard minimal element of ${\cal M}_n$.
Thus we have the following result of Lam and Leung.

\bbtw{\cite{Lam1}}{NonTw}
If $n\geqslant3$ is an integer, then
$\nu(n)=\xi(n)$ if and only if $n$ has at most two prime divisors.
\etw

Now, as a consequence of the previous considerations, we obtain:

\bwn{MinWnFin}
The number of a minimal set of generators of $k[X]^d$
is equal to $\xi(n)$
if and only if $n$ has at most two prime divisors.
\ewn

Note that in our examples all nonzero coefficients of the minimal
(standard or nonstandard) elements of ${\cal M}_n$ were equal to $1$.
Recently, John P. Steinberger \cite{Ste1} gave
the first explicit constructions
of nonstandard minimal elements of ${\cal M}_n$ (for some $n$)
with coefficients greater than $1$
(indeed containing arbitrary large coefficients).
He gave at the same time an answer to an old question of H.W. Lenstra Jr. \cite{Len}
concerning this subject.

\section{Polynomial constants of $\Delta$}\label{SecDelta}

\initpar
Let us recall that $\Delta$ is the derivation of $k[Y]$ given
by $\Delta(y_j)=y_j\left(y_{j+1}-y_j\right)$ for $j\in\mathbb Z_n$,
where $k[Y]=k[y_0,\dots,y_{n-1}]$.
It is a homogeneous derivation, that is, all the polynomials
$\Delta(y_0),\dots,\Delta(y_{n-1})$ are homogeneous of the same degree.
Put
$v=y_0y_1\cdots y_{n-1}$,
Observe that $v\in k[Y]^\Delta$. In this section we will prove
that $k[Y]^\Delta=k[v]$. For this aim  we first study Darboux polynomials of $\Delta$.

\smallskip

We say that a nonzero polynomial $F\in k[Y]$ is a {\it Darboux polynomial}
of $\Delta$, if $F$ is homogeneous and there exists a polynomial $\Lambda\in k[Y]$
such that $\Delta(F)=\Lambda F$. Such
a polynomial $\Lambda$ is uniquely determined and we say that $\Lambda$ is the {\it cofactor} of $F$.
Some basic properties of Darboux polynomials of arbitrary
homogeneous derivations one can find for example in
\cite{MoulinNS}, \cite{MN} or \cite{Now}.
Note that if $F,G\in k[Y]$ and $FG$ is a Darboux polynomial of $\Delta$,
then $F,G$ are also Darboux polynomials of $\Delta$ (\cite{MoulinNS}, \cite{Now}).
It is obvious that in our case each cofactor $\Lambda$ is of the form
$\lambda_0y_0+\lambda_1y_1+\dots+\lambda_{n-1}y_{n-1},$
where the coefficients $\lambda_0,\dots,\lambda_{n-1}$ belong to $k$.
We say that a Darboux polynomial is {\it strict} if it is not divisible by any of the variables
$y_0,\dots,y_{n-1}$.
The following  important proposition is a special case of Proposition $3$
from our paper \cite{MJN}.
For a sake of completeness we repeat its proof.

\bprop{DeltaDarb}
Let $F\in k[Y]\smallsetminus k$ be a strict Darboux polynomial of $\Delta$
and let $\Lambda=\lambda_0y_0+\dots+\lambda_{n-1}y_{n-1}$
be its cofactor. Then all $\lambda_i$ are integers
and they belong to the interval
$[-r,0]$, where $r=\deg F$. Moreover, two of the $\lambda_i$ at least are different from $0$.
\eprop

\begin{dowod}
As $F$ is strict, for any $i$, the polynomial $F_i=F_{|y_i=0}$
(that we get by evaluating $F$ in $y_i=0$)
is a nonzero homogeneous polynomial with the same degree $r$
in $n-1$ variables (all but $y_i$).
Evaluating the equality $\Delta(F)=\Lambda F$ at $y_{n-1}=0$ we obtain
$$
\sum_{i=0}^{n-3}y_i(y_{i+1}-y_i)\frac{\partial F_{n-1}}{\partial y_i}
-y_{n-2}^2 \frac{\partial F_{n-1}}{\partial y_{n-2}}=\left(\sum_{i=0}^{n-2}\lambda_iy_i\right)F_{n-1}.
\leqno{(\ast)}
$$
Let $r_0$ be the degree of $F_{n-1}$ with respect to $y_0$.
Then obviously $0\leqslant r_0\leqslant r$.
Consider now $F_{n-1}$ as a polynomial in $k[y_1,\dots,y_{n-2}][y_0]$.
Balancing monomials of degree $r_0+1$ in the equality $(\ast)$ gives
$\lambda_0=-r_0$. The same results hold for all coefficients of the cofactor $\Lambda$.

We already proved that all $\lambda_i$ are integers and $-r\leqslant \lambda_i\leqslant 0$.
Moreover, we proved that
$|\lambda_i|$ is the degree of $F_{i-1}$ with respect to $y_i$ (for any $i\in\mathbb Z_n$).
Thus $\lambda_i=0$ means that the variable $y_{i-1}$ appears in every monomial
of $F$ in which $y_i$ appears.
Then, if all $\lambda_i$ vanish, the product of all variables divides the nonzero
polynomial $F$,
a contradiction with the fact that $F$ is strict.
In the same way, if all $\lambda_i$ but one vanish, the variable corresponding to the
nonzero coefficient divides $F$, once again a contradiction.
\end{dowod}

\btw{TwDelta}
The ring of constants $k[Y]^\Delta$ is equal to $k[v]$, where $v=y_0y_1\dots,y_{n-1}$.
\etw

\begin{dowod}
The inclusion $k[v]\subseteq k[Y]^\Delta$ is obvious. We will prove the reverse inclusion.
For every Darboux polynomial $F$ of $\Delta$, we denote by $\Lambda(F)$
the cofactor of $F$. Then we have $\Delta(F)=\Lambda(F)\cdot F$, and
$\Lambda(F)=\lambda_0y_0+\dots+\lambda_{n-1}y_{n-1}$, where the coefficients
$\lambda_0,\dots,\lambda_{n-1}$ are uniquely determined. In  this case we denote
by $\Gamma(F)$ the sum $\lambda_0+\lambda_1+\dots+\lambda_{n-1}$.
In particular, the variables $y_0,\dots,y_{n-1}$ are Darboux polynomials of $\Delta$, and
$\Lambda(y_j)=y_{j+1}-y_j$, \  $\Gamma(y_j)=0$,
for any $j\in\mathbb Z_n$.
It follows from Proposition \ref{DeltaDarb} that
if a Darboux polynomial $F$ is strict and $F\not\in k$, then $\Gamma(F)$ is an integer, and
$\Gamma(F)\leqslant -2$. Note also that if $F,G$ are Darboux polynomials of $\Delta$, then
$FG$ is a Darboux polynomial of $\Delta$, and then
$$
\Lambda(FG)=\Lambda(F)+\Lambda(G)\quad\mbox{and}\quad
\Gamma(FG)=\Gamma(F)+\Gamma(G).
$$
Assume now that $F$ is a nonzero polynomial belonging to $k[Y]^\Delta$.
We will show that $F\in k[v]$. Since the derivation $\Delta$ is homogeneous
we may assume that $F$ is homogeneous.
Thus $F$ is a Darboux polynomial of $\Delta$ and its cofactor is equal to $0$.
Let us write this polynomial in the form
$$
F=y_0^{\beta_0}y_1^{\beta_1}\cdots y_{n-1}^{\beta_{n-1}}\cdot G,
$$
where $\beta_0,\dots,\beta_{n-1}$ are nonnegative integers, and $G$
is a nonzero from $K[Y]$ which is not divisible by any of the variables $y_0,\dots,y_{n-1}$.
Then $G$ is a strict Darboux polynomial of $\Delta$.
Let us suppose that $G\not\in k$.
Then $\Gamma(G)\leqslant -2$ (by Proposition \ref{DeltaDarb}), and we
have a contradiction:
$$
0=\Gamma(F)=\sum_{j=0}^{n-1}\beta_j\Gamma(y_j)+\Gamma(G)
=\sum_{j=0}^{n-1}\beta_j\cdot0+\Gamma(G)=\Gamma(G)\leqslant -2.
$$
Thus $F$ is a monomial of the form $b y^\beta=by_0^{\beta_0}y_1^{\beta_1}\cdots y_{n-1}^{\beta_{n-1}}$,
with some nonzero $b\in k$.
But $\Delta(F)=0$, so
$
\beta_0(y_1-y_0)+\beta_1(y_2-y_1)+\dots+\beta_{n-1}(y_0-y_{n-1})=0,
$
and so $\beta_0=\beta_1=\dots=\beta_{n-1}=c$, for some $c\in\mathbb N$.
Now we have  $F=by^\beta=b(y_0\cdots y_{n-1})^c=bv^c$,
and hence  $F\in k[v]$.
\end{dowod}

\section{The mappings @ and $\tau$}\label{SecAt}

\initpar
In this section we show that
the derivations $d$ and $\Delta$ have certain additional  properties,
and we present some specific relations between these derivations.

\smallskip

Let us fix the following two notations:
$$
\underline{a}=\left(\frac{x_1}{x_0},\frac{x_2}{x_1},\dots,\frac{x_{n-1}}{x_{n-2}}, \frac{x_0}{x_{n-1}}\right)
\quad\mbox{and}\quad v=y_0y_1\cdots y_{n-1}.
$$

We already know,
by Proposition \ref{HomPropA} and Theorem \ref{TwDelta},
that $k(X)^E=k\left(\underline{a}\right)$
and $k[Y]^\Delta=k[v]$.

\blem{ATlemA}
Let $F\in k[Y]$.
If $F(\underline{a})=0$, then there exists a polynomial $G\in k[Y]$ such that
$F=(v-1)G$.
\elem

\begin{dowod}
First note that if $b=(b_0,\dots,b_{n-1})$ is an element of $k^n$ such that
the product $b_0b_1\cdots b_{n-1}$ equals $1$, then $b$ is of the form
$
b=\left(\frac{c_1}{c_0},\frac{c_2}{c_1},\dots,\frac{c_{n-1}}{c_{n-2}}, \frac{c_0}{c_{n-1}}\right),
$
for some nonzero elements $c_0,\dots,c_{n-1}$ from $k$.
In fact, put: $c_0=1$, $c_1=b_0$, $c_2=b_0b_1,\dots,c_{n-1}=b_0b_1\cdots b_{n-2}$.

Let $P=v-1$, and let $A$ be the ideal of $\overline{k}[Y]=\overline{k}[y_0,\dots,y_{n-1}]$
generated by $P$, where $\overline{k}$ is the algebraic closure of $k$.
Observe that, for any $b\in \overline{k}^n$, if $P(b)=0$, then (by the assumption and the
above note) $F(b)=0$. This means, by the Nullstellensatz,
that some power of $F$ belongs to the ideal $A$.
But $A$ is a prime ideal, so $F\in A$ and so, there exists a polynomial $G\in\overline{k}[Y]$
such that $F=(v-1)G$. Since $F,\ v-1$ belong to $k[Y]$,
it is obvious that $G$ also belongs to $k[Y]$.
\end{dowod}

\blem{ATlemB}
Let $F$ is a nonzero homogeneous polynomial in $k[Y]$, then
$F(\underline{a})\neq0$.
\elem

\begin{dowod}
Suppose that $F(\underline{a})=0$.
Then, by Lemma \ref{ATlemA}, $F=(v-1)G$, for some $G\in k[Y]$.
As  $F$ is homogeneous,
the polynomials $v-1$ and $G$ are also homogeneous; but it is a contradiction,
because $v-1$ is not homogeneous.
\end{dowod}

Let us denote by $S$  the multiplicative subset $\left\{ F\in k[Y]; \  F(\underline{a})\neq 0\right\}$
and consider the quotient ring
$$
{\cal A}=S^{-1} k[Y].
$$
Every element of this ring is of the form $F/G$, where $F,G\in k[Y]$ and $G(\underline{a})\neq0$.
It is a local ring with the unique maximal ideal
$
I=\left\{\frac{F}{G}\in{\cal A}; \ F(\underline{a})=0\right\}.
$
It follows from Lemma~\ref{ATlemA} that $I=(v-1){\cal A}$.
Observe that $\Delta({\cal A})\subseteq {\cal A}$ and $\Delta(I)\subseteq I$, so
$\Delta$ is a derivation of ${\cal A}$ and $I$ is a differential ideal of ${\cal A}$.

If $f\in {\cal A}$, then $f(\underline{a})$ is well defined, and it is
a homogeneous rational function of degree zero,  that is, $f(\underline{a})\in k(X)^E$.
Thus we have a $k$-algebra homomorphism from ${\cal A}$ to $k(X)^E$.
This homomorphism we will denote by $@$. So we have:
$$
@:{\cal A}\to k(X)^E,\quad @(f)=f(\underline{a})\quad\mbox{for}\quad f\in{\cal A}.
$$
In particular, $@(v)=1$, and $@(y_j)=\frac{x_{j+1}}{x_j}$ for $j\in\mathbb Z_n$.
These equalities imply that $@$ is surjective. Note also that ker$@=I$, so
the field $k(X)^E$ is isomorphic to the factor ring~${\cal A}/I$.
Moreover, as a consequence of Lemma \ref{ATlemB} we have:

\bprop{ATlemProp}
If $f\in k(Y)$ is homogeneous and $@(f)=0$, then $f=0$.
\eprop

Note also the next important proposition.

\bprop{ATpropDD}
\ $\displaystyle d\circ@=@\circ \Delta$, that is,
$d\left(f(\underline{a})\right)=\left(\Delta(f)\right)(\underline{a})$
for   $f\in{\cal A}$.
\eprop

\begin{dowod}
It is enough to prove that the above equality holds
in the case when $f=y_j$ with $j\in\mathbb Z_n$.
Let $f=y_j$, $j\in\mathbb Z_n$. Then:
$$
\begin{array}{lcl}
d\Big(f(\underline{a})\Big)&=&d\left(\frac{x_{j+1}}{x_j}\right)
=\frac{d(x_{j+1})x_j-d(x_j)x_{j+1}}{x_j^2}
=\frac{x_{j+2}x_j-x_{j+1}^2}{x_j^2}
=\frac{x_{j+1}}{x_j}\left(\frac{x_{j+2}}{x_{j+1}}-\frac{x_{j+1}}{x_j}\right)\medskip\\
&=&\left(y_j(y_{j+1}-y_j)\right)(\underline{a})
=\left(\Delta(y_j)\right)(\underline{a})
=\left(\Delta(f)\right)(\underline{a}).
\end{array}
$$
This completes the proof.
\end{dowod}

\bwn{ATwn}
Let $f\in{\cal A}$.
If $\Delta(f)=0$, then $d\left(@(f)\right)=0$.
\ewn

\begin{dowod}
\ $d\left(@(f)\right)=@\left(\Delta(f)\right)=@(0)=0$ (by Proposition \ref{ATpropDD}).
\end{dowod}

Now we are ready to prove the following theorem.

\btw{ATtwPrime}
If $n$ is a prime number, then $k(Y)^\Delta=k(v)$, where
$v=y_0y_1\cdots y_{n-1}$.
\etw

\begin{dowod}
Put $P=v-1$. Note that $\Delta(P)=0$.
Let $0\neq f=\frac{F}{G}\in k(Y)$, where $F,G$ are nonzero, coprime
polynomials in $k[Y]$,
and assume that $\Delta(f)=0$. We will show, using an induction with respect to
$\deg F+\deg G$, that $f\in k(v)$.

If $\deg F+\deg G=0$, then $f\in k$, so $f\in k(v)$.
Assume that $\deg F+\deg G=r>0$.

If $P$ divides $F$, then $F=F'P$, for some $F'\in k[Y]$,
and then
$
\Delta\left(\frac{F'}{G}\right)=\frac{1}{P}\Delta\left(\frac{F}{G}\right)=0
$
with  $\deg F'+\deg G<r$.
Then, by induction, $\frac{F'}{G}\in k(v)$ and this implies
that $\frac{F}{G}\in k$, because
$\frac{F}{G}=P\frac{F'}{G}$ and $P\in k(v)$.
We use the same argument in the case when $P$ divides $G$.

Now we may assume that $P\nmid F$ and $P\nmid G$.
In this case, by Lemma \ref{ATlemA},
the quotient $\frac{F}{G}$ belongs to ${\cal A}$, and $@\left(\frac{F}{G}\right)\neq0$.
Moreover, we may assume that $\deg F\geqslant\deg G$
(in the opposite case we consider $G/F$ instead of $F/G$).

Since $\Delta(f)=0$, we have (by Corollary \ref{ATwn}) $@(f)\in k(X)^d\cap k(X)^E=k(X)^{d,E}$.
But $n$ is prime so, by Corollary \ref{HomWnPrime}, $k(X)^{d,E}=k$.
Therefore, $@\left(\frac{F}{G}\right)=c$, for some nonzero $c\in k$.
Thus we have
$$
\textstyle
0=@\left(\frac{F}{G}\right)-c=@\left(\frac{F}{G}-c\right)
=@\left(\frac{F-cG}{G}\right)=\frac{@(F-cG)}{@(G)},
$$
and hence, $@(F-cG)=0$. If $F-cG=0$, then $\frac{F}{G}=c\in k(v)$. Assume that $F-cG\neq 0$.
Then, by Lemma \ref{ATlemA}, $F-cG=H\cdot P$, for some nonzero $H\in k[Y]$.
As $\gcd(F,G)=1$, we have  $\gcd(H,G)=1$. Observe that $\Delta\left(\frac{H}{G}\right)=0$.
In fact,
$
\Delta\left(\frac{H}{G}\right)= \frac{1}{P}\Delta\left(\frac{PH}{G}\right)
= \frac{1}{P}\Delta\left(\frac{F-cG}{G}\right)
= \frac{1}{P}\Delta\left(\frac{F}{G}-c\right)
= \frac{1}{P}\Delta\left(\frac{F}{G}\right)=0.
$
It is clear that $\deg H+\deg G<\deg F+\deg G$.
Hence, by induction, the quotient $\frac{H}{G}$ belongs to $k(v)$.
But
$$
\textstyle
f=\frac{F}{G}=\left(\frac{F}{G}-c\right)+c=\frac{F-cG}{G}+c=P\frac{H}{G}+c,
$$
so $f\in k(v)$.
We proved that $k(Y)^\Delta\subseteq k(v)$.
The reverse inclusion is obvious.
\end{dowod}

Let us recall (see Theorem \ref{TwDelta}), that the ring of constants $k[Y]^\Delta$
is always equal to $k[v]$.
Thus, if $n$ is prime, then $k(Y)^\Delta$ is the field of quotients of $k[Y]^\Delta$.
In a general case a similar statement is not true. For example, if $n=4$,
then the rational function
$$
y_1y_3\frac{2y_0y_2-y_2y_3-y_0y_1}{y_1y_2+y_0y_3-2y_1y_3}
$$
belongs to $k(Y)^\Delta$ and it is not in $k(v)$.

\medskip

Let us recall (see Section \ref{SecNP})
that $\tau$ is an automorphism of $k(X)$ defined by
$$
\tau(x_j)=\E^jx_j
\quad
\mbox{for all}\quad j\in\mathbb Z_n.
$$

We say that a rational function $f\in k(X)$ is {\it $\tau$-homogeneous},
if $f$ is homogeneous in the ordinary sense and $\tau(f)=\E^s f$
for some $s\in\mathbb Z_n$.
In this case we say that $s$ is the $\tau$-degree of $f$ and we write
$\deg_{\tau}(f)=s$. Note that $\deg_{\tau}(f)$ is an element of $\mathbb Z_n$.

\smallskip

Let $\alpha=(\alpha_0,\dots,\alpha_{n-1})\in\mathbb Z^n$.
As usually, we denote by $x^\alpha$ the rational monomial
$x_0^{\alpha_0}\cdots x_{n-1}^{\alpha_{n-1}}$, and by $|\alpha|$
the sum $\alpha_0+\dots+\alpha_{n-1}$. Moreover, we denote by $\sigma(\alpha)$
the element from $\mathbb Z_n$ defined by
$$
\sigma(\alpha)=0 \alpha_0+1\alpha_1+2\alpha_2+\dots+(n-1)\alpha_{n-1} \ (\mbox{mod} \ n).
$$
Let us recall (see Section \ref{SecNP}) that
$\varrho:k(X)\to k(X)$ is a field  automorphism, defined by
$\varrho(x_j)=x_{j+1}$ for all $j\in\mathbb Z_n$.
It is very easy to check that:

\bprop{TAUrho}
Every rational monomial
$x^\alpha$, where $\alpha\in\mathbb Z^n$,
is $\tau$-homogeneous and its $\tau$-degree
is equal to $\sigma(\alpha)$. Moreover,
if $0\neq f\in k(X)$ and  $f$ is $\tau$-homogeneous,
then $\varrho(f)$ is also $\tau$-homogeneous,
and
 $
\MOD{\deg_{\tau} \varrho(f)}{\deg_\tau f+\deg f}{n}.
 $
\eprop

The derivation $d$ has the following additional properties.

\bprop{TAUpropA}
\ $\displaystyle \tau d \tau^{-1}=\E d$.
\eprop

\begin{dowod}
It is enough to show that $\tau d(x_j)=\E d(\tau(x_j))$ for $j\in\mathbb Z_n$.
Let us verify:
$
\tau d(x_j)=\tau (x_{j+1})=\E^{j+1}x_{j+1}=\E\cdot \E^j d(x_j)=\E d\left(\E^j x_j\right)=\E d(\tau(x_j)).
$
\end{dowod}

\bprop{TAUpropB}
Let $f\in k(X)$.
If $f$  is $\tau$-homogeneous, then $d(f)$ is $\tau$-homogeneous and
$
\deg_\tau d(f)=1+\deg_\tau f.
$
\eprop

\begin{dowod}
Assume that  $f$ is $\tau$-homogeneous and $s=\deg_\tau f$.
Since the derivation $d$ is homogeneous and $f$ is homogeneous in the ordinary sense,
$d(f)$ is also homogeneous in the ordinary sense.
Moreover, by the previous proposition, we have:
$\tau(d(f))=\E d(\tau(f))=\E d(\E^s f)=\E^{s+1} d(f)$,
so $d(f)$ is $\tau$-homogeneous and $\deg_\tau d(f)=s+1$.
\end{dowod}

\bprop{TAUpropC}
Let $F\in k[X]$ be a Darboux polynomial of $d$.
If $F$ is $\tau$-homogeneous, then $d(F)=0$.
\eprop

\begin{dowod}
Assume that $d(F)=bF$ with $b\in k[X]$, $F$ is homogeneous in the ordinary sense,
and $\tau(F)=\E^sF$.
Then $b\in k$, and  we have
$
\E d(F)=\E^{-s}\E d(\E^sF)=\E^{-s} \E d(\tau(F))=\E^{-s} \tau(d(F))
=\E^{-s} \tau(bF)=b\E^{-s} \tau(F) =b\E^{-s}\E^s F =bF=d(F).
$
Hence, $(\E-1)d(F)=0$. But $\E\neq1$, so $d(F)=0$.
\end{dowod}

\bprop{TAUpropD}
Let $f=\frac{P}{Q}$, where $P,Q$ are nonzero coprime
polynomials in $k[X]$.
If $f$ is $\tau$-homogeneous, then $P,Q$ are also $\tau$-homogeneous,
and
$\deg_\tau f=\deg_\tau P-\deg_\tau Q$.
Moreover, if $f$ is $\tau$-homogeneous and $d(f)=0$, then $d(P)=d(Q)=0$.
\eprop

\begin{dowod}
Assume that $f$ is $\tau$ homogeneous and $\deg_\tau f=s$.
Then $f$ is homogeneous in the ordinary sense and then, by Proposition \ref{HomHom},
the polynomials $P,Q$ are also homogeneous in the ordinary sense.
Since $\tau\left(\frac{P}{Q}\right)\E^s\frac{P}{Q}$, we have
$
\tau(P)Q=\E^sP\tau(Q)
$
and this implies that $\tau(P)=aP$, $\tau(Q)=bQ$, for some $a,b\in k[X]$
(because $P,Q$ are relatively prime). Comparing degrees, we deduce that $a,b\in k\smallsetminus\{0\}$.
But $\tau^n$ is the identity map, so $P=\tau^n(P)=a^nP$ and $Q=\tau^n(Q)=b^nQ$ and so,
$a$, $b$ are $n$-th roots of unity.
Since $\E$ is a primitive $n$-root, we have $a=\E^{s_1}$, $b=\E^{s_2}$,
for some $s_1,s_2\in\mathbb Z_n$. Thus, the polynomials $P,Q$ are $\tau$-homogeneous,
and it is clear that
$\MOD{s}{s_1-s_2}{n}$.

Assume now that $f$ is $\tau$-homogeneous and $d(f)=0$.
Then $P,Q$ are $\tau$-homogeneous Darboux polynomials of $d$
(with the same cofactor) and, by Proposition \ref{TAUpropC},
we have $d(P)=d(Q)=0$.
\end{dowod}

Note also the following proposition

\bprop{TAUpropE}
If $f\in k(Y)$ is  homogeneous, then
$@(f)$ is $\tau$-homogeneous, and
$\MOD{\deg_\tau @(f)}{\deg f}{n}$.
\eprop

\begin{dowod}
First assume that $f=F$ is a nonzero
homogeneous polynomial in $k[Y]$ of degree $s$
and consider all the monomial of $F$. Every nonzero monomial
is of the form $by^\alpha$, where $0\neq b\in k$, and $\alpha\in\mathbb N^n$ with $|\alpha|=s$.
For each such $y^\alpha$, we have  $@(y^\alpha)=x^\beta$, where
$
\beta=\Big(\beta_0,\dots,\beta_{n-1}\Big)=\Big(\alpha_{n-1}-\alpha_0, \ \alpha_0-\alpha_1, \ \alpha_1-\alpha_2, \ \dots, \ \alpha_{n-2}-\alpha_{n-1}\Big),
$
and then
$$\sigma(\beta)=\sum\limits_{j=0}^{n-1} j\beta_j
=|\alpha|-n\alpha_{n-1}=s - n\alpha_{n-1},$$
so $\MOD{\sigma(\beta)}{s}{n}$. This means that $\tau(x^\beta)=\E^s x^\beta$.
Thus, for  every nonzero monomial~$P$,
which appears in $F$, we have $\tau(@(P))=\E^s @(P)$.
This implies that $\tau(@(f))=\E^s @(f)$.
But $@(F)$ is also homogeneous
in the ordinary sense (because $@(F)\in k(X)^E$), so $@(F)$ is $\tau$-homogeneous,
and $\deg_\tau @(F)=\deg F$ (mod $n)$.

Now let $0\neq f\in k(Y)$ be an arbitrary  homogeneous rational function.
Let $f=\frac{F}{G}$ with $F,G\in k[Y]\smallsetminus\{0\}$ and
$\gcd(F,G)=1$.
Then $F,G$ are homogeneous
(by Proposition \ref{HomHom}), and
$@(f)=\frac{@(F)}{@(G)}$. Thus, by the above proof for polynomials,
$@(f)$ is $\tau$-homogeneous, and
$\MOD{\deg_\tau @(f)}{\deg f}{n}$.
\end{dowod}

\bprop{TAUpropF}
Let  $f,g\in k(Y)$ be homogeneous rational functions.
If $@(f)=@(g)$, then $f=v^c g$, for some $c\in\mathbb Z$.
\eprop

\begin{dowod}
Assume that $@(f)=@(g)$. Then, by Proposition \ref{TAUpropE},
$
\MOD{\deg f}{\deg_\tau @(f)=\deg_\tau @(g)\equiv \deg g}{n},
$
so there exists $c\in\mathbb Z$ such that $\deg f=nc+\deg g$.
Then $f$ and $v^cg$ are homogeneous of the same degree, so $f-v^cg$ is homogeneous.
Observe that $@(f-v^cg)=@(f)-@(v)^c@(g)=@(f)=@(g)=0$.
Hence, by Proposition \ref{ATlemProp}, we have $f=v^cg$.
\end{dowod}

\medskip

Let us assume that $g$ is a $\tau$-homogeneous rational function belonging
to the field $k(X)^{d,E}$. We will show that then there exists
a homogeneous (in the ordinary sense) rational function $f\in k(Y)$ such
that $\Delta(f)=0$ and $@(f)=g$. This fact will play a key role
in our description of the structure of the field $k(Y)^\Delta$.
For a proof of this fact we need to prove some
lemmas and propositions

\medskip

Let us recall from Section \ref{SecNP}, that the elements $e_0,\dots,e_{n-1}\in\mathbb Z^n$ are defined by:
$e_0=(1,0,0,\dots,0)$, $e_1=(0,1,0,\dots,0), \ \dots, \ e_{n-1}=(0,0,\dots,0,1)$.
In particular, we have
$$
@(y_j)=\frac{x_{j+1}}{x_j}=x^{e_{j+1}-e_j},\quad\mbox{for}\quad j\in\mathbb Z_n.
$$

\blem{TAUlemEE}
Let $\alpha\in\mathbb Z^n$. Assume that $|\alpha|=0$ and $\sigma(\alpha)=0$ $(\mbox{{\em mod}} \ n)$.
Then there exist a sequence $\beta=(\beta_0,\dots,\beta_{n-1})\in\mathbb Z^n$
such that $|\beta|=0$ and
$
 \alpha=\sum_{j=0}^{n-1}\beta_j(e_{j+1}-e_j).
$
\elem

\begin{dowod}
Since $\MOD{\sigma(\alpha)}{0}{n}$, there exists an integer $r$ such that
$n\alpha_0+\sigma(\alpha)=-rn$.
Put: $\beta_0=r$ and
$
\beta_j=r-\sum_{i=1}^j\alpha_i,
$
for $j=1,\dots,n-1$.
\end{dowod}

\blem{TAUlemXY}
If  $\alpha\in\mathbb Z^n$ with $|\alpha|=0$, then there exists $\beta\in\mathbb Z^n$
such that $@(y^\beta)=x^\alpha$.
\elem

\begin{dowod}
Put: $\beta_j=\sum_{i=j+1}^{n-2} \alpha_i$ \  for $j=0,1,\dots,n-3$,
 \ and $\beta_{n-2}=0$, \ $\beta_{n-1}=-\alpha_{n-1}$.
\end{dowod}

\medskip

Now we assume that $P$ is a fixed nonzero
$\tau$-homogeneous polynomial in $k[X]$.
Let us write this polynomial in the form
$$
P=c_1 x^{\gamma_1}+\dots+c_r x^{\gamma_r},
$$
where $c_1,\dots,c_r$ are nonzero elements of $k$, and $\gamma_1,\dots,\gamma_r\in\mathbb N^n$.
For every $q\in\{1,\dots,r\}$,
we have
$|\gamma_q|=\deg F$
and $\MOD{\sigma(\gamma_q)}{\deg_\tau F}{n}$, and hence,
$|\gamma_q-\gamma_1|=0$ \ and  $\MOD{\sigma(\gamma_q-\gamma_1)}{0}{n}$.
This implies, by Lemma \ref{TAUlemEE}, that for any $q\in\{1,\dots,r\}$,
there exists a sequence $\beta^{(q)}=\left(\beta_0^{(q)},\dots,\beta_{n-1}^{(q)}\right)\in\mathbb Z^n$
such that $\left|\beta^{(q)}\right|=0$ and
$$
\gamma_q-\gamma_1=\sum_{j=0}^{n-1}\beta_j^{(q)}\left(e_{j+1}-e_j\right).
$$
For each $j\in\{0,1,\dots,n-1\}$, we define:
$$
\alpha_j=\min\left\{\beta_j^{(1)}, \ \beta_j^{(2)}, \ \dots, \beta_j^{(r)}\right\},
$$
and we denote by $\lambda$ the sequence $(\lambda_0,\dots,\lambda_{n-1})\in\mathbb Z^n$ defined by
$$
\lambda=\gamma_1+\sum_{j=0}^{n-1}\alpha_j\left(e_{j+1}-e_j\right).
$$
Observe that $|\lambda|=|\gamma_1|=\deg P$, and
$\gamma_q=\lambda+\sum_{j=0}^{n-1}\left(\beta_j^{(q)}-\alpha_j\right)\left(e_{j+1}-e_j\right)$
for any $q\in\{1,\dots,r\}$, and moreover, each
$\beta_j^{(q)}-\alpha_j$ is a nonnegative integer.
Put
$a_{qj}=\beta_j^{(q)}-\alpha_j$, for $j\in\mathbb Z_n$, $q\in\{1,\dots,r\}$,
and
$a_q=\left(a_{q0}, \ a_{q1}, \ \dots, \ a_{q(n-1)}\right)$
for all $q=1,\dots,r$.
Then each $a_q$ belongs to $\mathbb N^n$, and we have the equalities
$$
\gamma_q=\lambda+\sum_{j=0}^{n-1} a_{qj}\left(e_{j+1}-e_j\right),
\quad\mbox{for any}\quad q\in\{1,\dots,r\}.
$$
Let us remark that $\lambda \in \mathbb N^n$.

\noindent Indeed, for any $j\in \mathbb Z_n$,
we have $\lambda_j = \gamma_{1j} + \alpha_{j-1} - \alpha_j$, where $\alpha_{j-1}=\beta^{(q)}_{j-1}$
for some $q$ and $\alpha_j \leq \beta^{(q)}_j$.
Thus $\lambda_j = \gamma_{1j} + \beta^{(q)}_{j-1} - \alpha_j \geq
\lambda_j = \gamma_{1j} + \beta^{(q)}_{j-1} - \beta^{(q)}_j = \gamma_{qj} \geq 0$.\smallskip

\noindent Moreover, $|a_q|=|\beta^{(q)}-\alpha|=|\beta^{(q)}|-|\alpha|=-|\alpha|$,
because $|\beta^{(q)}|=0$. This means that $|\alpha|\leqslant 0$, and all the
numbers $|a_1|, \dots,|a_r|$ are the same;  they are equal to $-|\alpha|$.
Consider the polynomial in $k[Y]$ defined by
$$
\overline{P}=c_1y^{a_1}+\cdots+c_ry^{a_r}.
$$
It is a nonzero homogeneous (in the ordinary sense) polynomial
of degree $-|\alpha|$. It is easy to check that
$@(\overline{P})=x^{-\lambda}P$. Thus, we proved the following proposition.

\bprop{TAUpropLem}
If $P\in k[X]$ is a nonzero $\tau$-homogeneous polynomial, then there
exist a sequence $\lambda\in\mathbb Z^n$ and a homogeneous polynomial $\overline{P}\in k[Y]$
such that $@(\overline{P})=x^{-\lambda}P$ and $|\lambda|=\deg P$.
\eprop

\buw{TAUuw}
\begin{em}
In the above construction, the polynomial $\overline{P}$ is not divisible by
any of the variables $y_0,\dots,y_n$.
Let us additionally assume that $d(P)=0$. Then it is not difficult to show
that
$
\Delta(\overline{P})=-(\lambda_0y_0+\dots+\lambda_{n-1}y_{n-1}) \overline{P},
$
that is, $\overline{P}$ is a strict Darboux polynomial of $\Delta$ and its
cofactor is equal to $-\sum \lambda_i y_i$.
This implies, by Proposition \ref{DeltaDarb}, that if additionally $d(P)=0$,
among all nonnegative numbers $\lambda_0,\dots,\lambda_{n-1}$, at least two are
different from zero.
\end{em}
\euw

Now we are ready to prove the following, mentioned above, proposition.

\bprop{TAUpropGlowne}
Let  $g$ be a $\tau$-homogeneous rational function belonging
to the field $k(X)^{d,E}$.
Then there exists
a homogeneous rational function $f\in k(Y)$ such
that $\Delta(f)=0$ and $@(f)=g$.
\eprop

\begin{dowod}
For  $g=0$ it is obvious. Assume that $g\neq0$, and let $g=\frac{P}{Q}$,
where $P,Q\in k[X]\smallsetminus\{0\}$ with $\gcd(P,Q)=1$.
It follows from Propositions \ref{HomHom} and \ref{TAUpropD}, that
the polynomials $P,Q$ are homogeneous (in the ordinary sense)
of the same degree, and they are also $\tau$-homogeneous.
By Proposition \ref{TAUpropLem},
there
exist  sequences $\lambda, \mu \in\mathbb Z^n$ and a
homogeneous polynomials $\overline{P}, \overline{Q}\in k[Y]$
such that
$@(\overline{P})=x^{-\lambda}P$, $@(\overline{Q})=x^{-\mu}Q$, and
$|\lambda|=|\mu|=\deg P=\deg Q$. Then we have
$$
g=\frac{P}{Q}=\frac{x^\lambda\left(x^{-\lambda}P\right)}{x^\mu\left(x^{-\mu}Q\right)}
=\frac{x^\lambda @(\overline{P})}{x^\mu @(\overline{Q})}
=x^{\lambda-\mu}@\left(\overline{P}/\overline{Q}\right).
$$
Since $|\lambda-\mu|=0$,
there exists (by Lemma \ref{TAUlemXY})
$\beta\in\mathbb Z^n$ such that $@(y^\beta)=x^{\lambda-\mu}$.
Put
$
f=y^\beta\cdot \overline{P}/\overline{Q}.
$
Then $f\in k(Y)$ is a homogeneous rational function, and $@(f)=g$.
Now we will show that $\Delta(f)=0$.
To this aim let us recall that $g$ belongs to the field $k(X)^{d,E}$,
so $d(g)=0$. This implies  that
$@(\Delta(f))=0$, because (by Proposition \ref{ATpropDD})
$@(\Delta(f))=d(@(f))=d(g)=0$.
But the rational function $\Delta(f)$ is  homogeneous, so
by Proposition \ref{ATlemProp}, $\Delta(f)=0$.
\end{dowod}

\section{Rational constants of $\Delta$}\label{SecRC}

\initpar
We proved (see Proposition \ref{HomDerZero}) that $k(X)^{d,E}=k(q_1,\dots,g_{m-1})$,
where $m=n-\f(n)$, and $g_1,\dots,g_{m-1}\in k(X)$
are some algebraically independent
homogeneous rational functions of degree $0$.
We proved in  fact, that each
$g_j=$ (for $j=1,\dots,m-1$)
is equal to the quotient $\frac{w_j}{w_0}$.
These quotients are usually not $\tau$-homogeneous.
We will show in the next section
that,  in some cases, we are ready to find such algebraically independent
generators of $k(X)^{d,E}$
which are  additionally $\tau$-homogeneous.
In this section we prove that if we have $\tau$-homogeneous generators, then
we may construct some algebraically independent generators of the field $k(Y)^\Delta$.

\medskip

Let us assume that $k(X)^{d,E}=k(g_1,\dots,g_{m-1})$, where $g_1,\dots,g_{m-1}\in k(X)$
are algebraically independent $\tau$-homogeneous rational functions.
We know, by Proposition \ref{TAUpropGlowne}, that
for each $g_j$ there exists a homogeneous rational function $f_j\in k(Y)$ such that
$\Delta(f_j)=0$ and $@(f_j)=g_j$. Thus we have homogeneous rational
functions $f_1,\dots,f_{m-1}$, belonging to the field $k(Y)^\Delta$.
We know also that $v\in k(Y)^\Delta$, where $v=y_0y_1\cdots y_{n-1}$.
In this section we will prove the following theorem.

\btw{RCtw}
Let $g_1,\dots,g_{m-1}$ and $v, f_1,\dots,f_{m-1}$ be as above.
Then the elements $v,f_1,\dots,f_{m-1}$
are algebraically independent over $k$, and
$
k(Y)^\Delta=k(v,f_1,\dots,f_{m-1}).
$
\etw

We will prove it in several steps.

\medskip

{\bf Step 1.}
{\it The elements $f_1,\dots,f_{m-1}$ are algebraically independent over $k$}.

\medskip

\begin{dowod}
Suppose that $W(f_1,\dots,f_{m-1})=0$ for some $W\in k[t_1,\dots,t_{m-1}]$.
Then
$$
0=@\Big(W(f_1,\dots,f_{m_1})\Big)=W\Big(@(f_1),\dots,@(f_{m-1})\Big)=
W(g_1,\dots,g_{m-1}).
$$
But $g_1,\dots,g_{m-1}$ are algebraically independent, so $W=0$.
\end{dowod}

In the next steps we  write $f$ instead of $\{f_1,\dots,f_{m-1}\}$,
and $g$ instead of $\{g_1,\dots,g_{m-1}\}$
In particular, $k(f)$ means $k(f_1,\dots,f_{m-1})$,

\medskip

{\bf Step 2.}
{\it $v\not\in k(f)$}.

\medskip

\begin{dowod}
Suppose that $v\in k(f)$.  Let $v=P(f)/Q(f)$ for some
$P,Q\in k[t_1,\dots,t_{m-1}]$. Then $Q(f)v-P(f)=0$ and we have
$0=@(Q(f)v-P(f))=Q(g)@(v)-P(g)$. But $@(v)=1$, so $P(g)=Q(g)$, and so $P=Q$,
because $g_1,\dots,g_{m-1}$ are algebraically independent.
Thus $v=P(f)/Q(f)=P(f)/P(f)=1$; a contradiction.
\end{dowod}

{\bf Step 3.}
{\it The elements $v, f_1,\dots,f_{m-1}$ are algebraically independent over $k$}.

\medskip

\begin{dowod}
We already know (by Step $1$) that $f_1,\dots,f_{m-1}$ are algebraically independent.
Suppose that $v$ is algebraic over $k(f)$.
Let $F(t)=b_rt^r+\dots+b_1t+b_0\in k(f)[t]$  (with $a_r\neq0)$
be the minimal polynomial of $v$ over $k(f)$.
Multiplying by the common denominator, we may assume that the coefficients
$b_0,\dots,b_r$ belong to the ring $k[f]$.
There exists polynomials $B_0,B_1,\dots,B_r\in k[t_1,\dots,t_{m-1}]$
such that $b_j=B_j(f)$ for all $j=0,\dots,r$.
Thus,
$B_r(f)v^r+\dots+B_1(f)v+B_0(f)=0.$
Using $@$, we obtain the equality
$$
B_r(g)1^r+\dots+B_1(g)1+B_0(g)=0,
$$
which implies that
$
B_r+\cdots+B_1+B_0=0,
$
because $g_1,\dots,g_{m-1}$ are algebraically independent over $g$.
This means, in particular, that $F(1)=0$.
But $F(t)$ is an irreducible polynomial of degree $r\geqslant1$,
so $r=1$. Hence, $B_1(f)v+B_0(f)=0$, $B_1(f)\neq0$,
and hence $v=-B_0(f)/B_1(f)\in k(f)$; a contradiction with Step $2$.
\end{dowod}

It is clear that $k(v,f)\subseteq k(Y)^\Delta$. For a proof of Theorem \ref{RCtw}
we must show that the reverse inclusion also holds.
Note that the derivation $\Delta$ is homogeneous,
so it is well known that its
field of constants is generated by some
homogeneous rational functions.
Hence for a proof of this theorem we need to prove that every homogeneous
element of $k(Y)^\Delta$ is an element of $k(v,f)=k(v,f_1,\dots,f_{m-1})$.

\medskip

Let us assume that $H$ is a nonzero
homogeneous rational function belonging to $k(Y)^\Delta$,
and put $h=@(H)$.

\medskip

{\bf Step 4.}
{\it $h\in k(g)$ and $h$ is $\tau$-homogeneous}.

\medskip

\begin{dowod}
Since $h=@(H)$, we have $h\in k(X)^E$. Moreover, $d(h)=d@(H)=@\Delta(H)=@(0)=0$, so
$h\in k(X)^d\cap k(X)^E=k(X)^{d,E}=k(g)$.
The $\tau$-homogeneity of $h$ follows from Proposition \ref{TAUpropE}.
\end{dowod}

Now we introduce some new notations.
The $\tau$-degrees of $g_1,\dots,g_{m-1}$ we denote by $s_1,\dots,s_{m-1}$,
respectively, and by $s$ we denote the $\tau$-degree of $h$.
Thus we have $\tau(g_j)=\E^{s_j} g_j$ for $j=1,\dots,m-1$, and $\tau(h)=\E^s h$.
We already know that $h\in k(g)$, so we have
$$
h=\frac{A(g)}{B(g)}
$$
for some relatively prime nonzero
polynomials $A,B\in k[t_1,\dots,t_{m-1}]$.

\medskip

{\bf Step 5.}
{\it The elements $A(g), \ B(g)$ are $\tau$-homogeneous}.

\medskip

\begin{dowod}
Since $\tau(h)=\E^s h$, we have
$\tau(A(g)) B(g)=\E^s A(g)\tau(B(g))$, that is,
$$
A\Big(\E^{s_1}g_1,\dots,\E^{s_{m-1}} g_{m-1}\Big) B\Big(g_1,\dots,g_{m-1}\Big)=
\E^s A\Big(g_1,\dots,g_{m-1}\Big) B\Big(\E^{s_1}g_1,\dots,\E^{s_{m-1}} g_{m-1}\Big).
$$
But the elements $g_1,\dots,g_{m-1}$ are algebraically independent over $k$, so in the polynomial ring
$k[t_1,\dots,t_{m-1}]$ we have the equality
$$
A\Big(\E^{s_1}t_1,\dots,\E^{s_{m-1}} t_{m-1}\Big)\cdot B=
\E^s A\cdot  B\Big(\E^{s_1}t_1,\dots,\E^{s_{m-1}} t_{m-1}\Big),
$$
which implies that $A\Big(\E^{s_1}t_1,\dots,\E^{s_{m-1}} t_{m-1}\Big)=pA$
and $B\Big(\E^{s_1}t_1,\dots,\E^{s_{m-1}} t_{m-1}\Big)=qB$, for some
$p,q\in k[t_1,\dots,t_{m-1}]$ (because we assumed that $\gcd(A,B)=1$).
Comparing degrees we deduce that $p,q\in k$.
Therefore,
$\tau(A(g))=A(\tau(g_1,\dots,\tau(g_{m-1}))=
A\Big(\E^{s_1}g_1,\dots,\E^{s_{m-1}} g_{m-1}\Big) =pA(g_1,\dots,g_{m-1})=pA(g)$,
so, $\tau(A(g))=p A(g)$, and similarly $\tau(B(g))=qB(g)$.
But $\tau^n$ is the identity map, so $p^n=q^n=1$ and so, $p,q$ are $n$-th roots of unity.
Put  $p=\E^a$ and $q=\E^b$, where $a,b\in\mathbb Z_n$.
Then we have $\tau(A(g))=\E^a A(g)$ and $\tau(B(g))=\E^b B(g)$.
Moreover, $A(g)$, $B(g)$ are homogeneous in the ordinary sense, because
they belong to $k(X)^E$, so they
are homogeneous rational functions of degree zero.
This means that $A(g)$, $B(g)$ are $\tau$-homogeneous.
\end{dowod}

Let us fix: $a=\deg_\tau A(g)$ and $b=\deg_\tau B(g)$.

\medskip

If $\alpha=(\alpha_1,\dots,\alpha_{m-1})\in\mathbb N^{m-1}$ then, as usually,
we denote by $t^\alpha$ and $g^\alpha$
the elements $t_1^{\alpha_1}\cdots t_{m-1}^{\alpha_{m-1}}$ and
$g_1^{\alpha_1}\cdots g_{m-1}^{\alpha_{m-1}}$, respectively,
and moreover, we denote:
$$
\begin{array}{lcl}
w(\alpha)&=&\alpha_1 s_1+\dots+\alpha_{m-1} s_{m-1},\smallskip\\
u(\alpha)&=&\alpha_1 \deg f_1+\dots+\alpha_{m-1} \deg f_{m-1}.
\end{array}
$$
Recall that $s_j=\deg_\tau(g_j)$ and $@(f_j)=g_j$,  for all $j=1,\dots, m-1$.
It follows from Proposition \ref{TAUpropE} that for each $j$ we have the congruence
$\MOD{s_j}{\deg f_j}{n}$.
Therefore,
$$
\MOD{u(\alpha)}{w(\alpha)}{n}\quad\mbox{for all}\quad \alpha\in\mathbb N^{n-1}.
$$

Let us write the polynomials $A,B$ in the forms
$$
A=\sum_{\alpha \in S_A} A_\alpha t^\alpha,\quad
B=\sum_{\beta \in S_B} B_\beta t^\beta,
$$
where $A_\alpha$, $B_\beta$ are nonzero elements of $k$, and $S_A$, $S_B$
are  finite  subsets of $\mathbb N^{m-1}$.

\medskip

{\bf Step 6.}
{\it
 \ $\MOD{w(\alpha)}{a}{n}$ for all $\alpha\in S_A$,
and $\MOD{w(\beta)}{b}{n}$ for all $\beta\in S_B$
}.

\medskip

\begin{dowod}
Since $\tau(A(g))=\E^a A(g)$, we have
$$
\begin{array}{lcl}
\E^a \sum A_\alpha g^\alpha&=&\E^aA(g)=\tau(A(g))=\sum A_\alpha \tau\left(t^\alpha\right)\medskip\\
&=&\sum A_\alpha \left(\E^{s_1} g_1\right)^{\alpha_1}\cdots\left(\E^{s_{m-1}} g_{m-1}\right)^{\alpha_{m-1}}\medskip\\
&=&\sum A_\alpha \E^{w(\alpha)}g^\alpha.

\end{array}
$$
Hence, $\sum A_\alpha\left(\E^a-\E^{w(\alpha)}\right) g^\alpha=0$.
But $g_1,\dots,g_{m-1}$ are algebraically independent and each $A_\alpha$ is nonzero,
so $\E^{w(\alpha)}=\E^a$ and consequently $\MOD{w(\alpha)}{a}{n}$,
for all $\alpha\in S_A$.
The same we do for the elements $w(\beta)$.
\end{dowod}

Since $\MOD{u(\alpha)}{w(\alpha)}$ for all $\alpha\in\mathbb N^{m-1}$,
it follows from the above step that, for each $\alpha\in S_A$, there exists
$p(\alpha)\in \mathbb Z$ such that $u(\alpha)=a+p(\alpha)n$. Put
$$
p=\max\left(\{0\}\cup \{p(\alpha); \ \alpha\in S_A\}\right),
$$
and put $a(\alpha)=p-p(\alpha)$ for $\alpha\in S_A$.
Then all $a(\alpha)$ are nonnegative integers
and all the numbers  $u(\alpha)+a(\alpha)n$, for each $\alpha\in S_A$,
are the same; they are equal to $a+pn$.

A similar procedure we do with elements of  $S_B$.
For each $\beta\in S_B$ there exists an integer $b(\beta)$ such that
$
u(\beta)+b(\beta)n=b+qn,
$
for all $\beta\in S_B$, where $q$ is a nonnegative integer.
Consider now the following quotient
$$
\Theta=\frac{\sum\limits_{\alpha\in S_A} A_\alpha f^\alpha v^{a(\alpha)}}
{\sum\limits_{\beta\in S_B} B_\beta f^\beta v^{b(\beta)}}.
$$
This quotient belongs of course to  $k(v,f_1,\dots,f_{n-1})$. In its
numerator each component $A_\alpha f^\alpha v^{a(\alpha)}$, for all $\alpha\in S_A$,
is a homogeneous
rational function of the same degree $a+pn$, so  the numerator is homogeneous.
By the same way we see that the denominator is also homogeneous.
Hence, $\Theta$ is a homogeneous rational function.
Observe that $@(\Theta)=h$.
We have also $@(H) =h$. Thus,
$H$ and  $\Theta$ are two homogeneous rational functions
such that  $@(H)=@(\Theta)$. By Proposition \ref{TAUpropF}, there exists an integer
$c$ such that
$
H=v^c\cdot\Theta.
$
Therefore, $H\in k(v,f_1,\dots,f_{n-1})$.
This completes our proof of Theorem \ref{RCtw}. $\square$

\section{Two special cases}\label{SecRPS}

\initpar
In this section we present a description of the field
$k(Y)^\Delta$ in the case when $n$ is a power of a prime number,
and in the case when $n$ is a product of two primes.

\smallskip

Let $n=p^s$, where $p$ is   prime  and $s\geqslant1$.
We already know, by Theorem \ref{ATtwPrime}, that if $s=1$, then
$k(Y)^\Delta=k(v)$. Now we assume that $s\geqslant2$.

\btw{RPStw}
If $n=p^s$, where $p$ is prime and $s\geqslant2$, then
$$
k(Y)^\Delta=k(v, f_1,\dots,f_{m-1})
$$
with $m=p^{s-1}$, where $v=y_0\cdots y_{n-1}$ and $f_1,\dots,f_{m-1}\in k(Y)$
are homogeneous rational functions such that
$v,f_1,\dots,f_{m-1}$ are algebraically independent over $k$.
\etw

\begin{dowod}
In this  case $m=n-\f(n)=p^s-\f(p^s)=p^{s-1}$ and hence,  $n=pm$.
Since
$
\Phi_{p^s}(t)=1+t^m+t^2m+\dots+t^{(p-1)m},
$
we have:
$
w_0=u_0 u_m u_{2m}\cdots u_{(p-1)m},
$
and
$w_j=u_{0m+j}u_{1m+j}u_{2m+j}\cdots u_{(p-1)m+j},
$
for all $j=0,1,\dots,m-1$.
Recall (see Lemma \ref{NVlemC}) that $\tau(u_j)=u_{j+1}$
for $j\in\mathbb Z_n$, so each
$w_j$ is equal to $\tau^j(w_0)$.

\smallskip

Observe that
$
\tau^m(w_0)=w_0.
$
This implies that the $\tau$-degree of
every nonzero monomial (with respect to variables $x_0,\dots,x_{n-1})$
of $w_0$ is divisible by $p$.
This means that in the $\tau$-decomposition of $w_0$ there are only
components with $\tau$-degrees $0,p,2p,\dots,(m-1)p$.
Let
$
w_0=v_0+v_1+\dots+v_{m-1},
$
where each $v_j\in k[X]$ is  $\tau$-homogeneous and
$\tau(v_j)=\E^{pj}v_j$.
Of course $d(v_j)=0$ for all $j$ (because $\tau d=\E d\tau$),
and  $\deg (v_j)=p$ for all $j$ (by Proposition~\ref{DdHom}).
Now observe
that
if $p\geqslant3$ then $\varrho(w_0)=w_0$, and if $p=2$ then
$\varrho(w_0)=-w_0$. Hence
$\varrho(w_0)=\pm w_0$, and we have
$$
v_0+v_1+\dots+v_{m-1}=w_0=\pm \varrho(w_0)=\pm(\varrho(v_0)\pm\varrho(v_1)\pm\dots\pm
\varrho(v_{m-1})
$$
Since the $\tau$-decomposition of $w_0$ is unique,
we deduce (by Proposition \ref{TAUrho}), that
$$
v_1=\pm\varrho(v_0),\quad\ v_2=\pm\varrho(v_1),\quad \dots,\quad  v_{m-1}=\pm\varrho(v_{m-2}),\quad v_0=\pm\varrho(v_{m-1}),
$$
and we have  $v_j=\pm \varrho^j(v_0)$ for all $j=0,1,\dots,m-1$.
Therefore, the $\tau$-decomposition of $w_0$ is of the form
$
w_0=v_0+b_1\varrho(v_0)+b_2\varrho^2(v_0)+\dots+b_{m-1}\varrho^{m-1}(v_0),
$
where the coefficients $b_1,\dots,b_{m-1}$ belong to $\{-1,1\}$. This implies that
$$
w_1=\tau(w_0)=v_0+b_1\E^p\varrho(v_0)+b_2\E^{2p}\varrho^2(v_0)+\dots+b_{m-1} \E^{(m-1)p}\varrho^{m-1}(v_0).
$$
We do the same for $w_2=\tau(w_1)=\tau^2(w_0)$, and for all $w_j$.
Thus,
for all $j=0,1,\dots,m-1$, we have
$
w_j=v_0+c_{j1}\varrho(v_0)+c_{j2}\varrho^2(v_0)+\dots+c_{j,m-1} \varrho^{m-1}(v_0),
$
where each $c_{ji}$ belongs to the ring $\mathbb Z[\E]$.
Consider now the rational functions $g_1,\dots,g_{m-1}\in k(X)$ defined by
$$
g_j=\frac{\varrho^j(v_0)}{v_0},
$$
for $j=1,\dots,m-1$.
These functions are $\tau$-homogeneous. They are homogeneous of degree zero,
and they are constants of $d$.
Moreover, if $j\in\{1,\dots,m-1\}$, then we have:
$$
\frac{w_j}{w_0}=\frac{v_0+\sum\limits_{i=1}^{m-1}c_{ji}\varrho^i(v_0)}
{v_0+\sum\limits_{i=1}^{m-1}c_{0i}\varrho^i(v_0)}
=\frac{1+v_0^{-1}\sum\limits_{i=1}^{m-1}c_{ji}\varrho^i(v_0)}
{1+v_0^{-1}\sum\limits_{i=1}^{m-1}c_{0i}\varrho^i(v_0)}
=
\frac{1+\sum\limits_{i=1}^{m-1}c_{ji}g_i}
{1+\sum\limits_{i=1}^{m-1}c_{0i}g_i}.
$$
Hence, all the elements $\frac{w_1}{w_0}, \dots,\frac{w_{m-1}}{w_0}$
belong to the field $k(g_1,\dots,g_{m-1})$,
and hence, by Proposition \ref{HomDerZero}, the elements $g_1,\dots,g_{m-1}$ are algebraically
independent over $k$ and we have the equality
$
k(X)^{E,d}=k(g_1,\dots,g_{m-1}).
$
Note that $g_1,\dots,g_{m-1}$ are $\tau$-homogeneous.
It follows from Proposition \ref{TAUpropGlowne}, that
for each $g_j$ there exists
a homogeneous rational function $f_j\in k(Y)$ such that
$\Delta(f_j)=0$ and $@(f_j)=g_j$.
We know, by Theorem \ref{RCtw}, that the elements
$v, f_1,\dots,f_{m-1}$, are algebraically independent over $k$,
and
$
k(Y)^\Delta=k(v,f_1,\dots,f_{m-1}).
$
This completes our proof of Theorem \ref{RPStw}.
\end{dowod}

Using the above theorem and its proof we obtain:

\bexa{RPSexa}
If $n=4$, then $k(Y)^\Delta=k(v,f)$, where
$f=y_1y_3\frac{2y_0y_2-y_2y_3-y_0y_1}{y_1y_2+y_0y_3-2y_1y_3}$
and $v=y_0y_1y_2y_3$.
\eexa


Consider the case $n=6$.

\bexa{PQexa}
If $n=6$, then $k(Y)^\Delta=k(v,f_1,f_2,f_3)$,
where $v=y_0\cdots y_{5}$, and $f_1,f_2,f_3$
are some homogeneous rational functions in $k(Y)$ such that
$v,f_1, f_2,f_{3}$ are algebraically independent over $k$.
\eexa

\begin{dowod}
We have:
$\f(n)=\f(6)=2$, $m=n-\f(n)=4$,
$\Phi_6(t)=t^2-t+1$, and
$w_0=\frac{u_0u_2}{u_1}$, \  $w_1=\frac{u_1u_3}{u_2}=\tau(w_0)$,
\ $w_2=\frac{u_2u_4}{u_3}=\tau^2(w_0)$, \ $w_3=\frac{u_3u_5}{u_4}=\tau^3(w_0)$.
Let us denote:
$F_0=u_0u_2u_4, \quad F_1=u_1u_3u_5=\tau(F_0),\quad
G_0=u_0u_3,\quad G_1=u_1u_4=\tau(G_0),\quad G_2=u_2u_5=\tau^2(G_0).
$
It is clear that the polynomials $F_0,F_1,G_0,G_1,G_2$
are constants of  $d$.
Note that $w_0=\frac{F_0}{G_1}$, \ $w_1=\frac{F_1}{G_2}$,
\ $w_2=\frac{F_0}{G_0}$, \ $w_3=\frac{F_1}{G_1}$,
so we have:
$
\frac{w_1}{w_0}=\frac{F_1 G_1}{F_0G_2},\quad
\frac{w_2}{w_0}=\frac{F_0G_1}{F_0G_0}=\frac{G_1}{G_0},\quad
\frac{w_3}{w_0}=\frac{F_1G_1}{F_0G_1}=\frac{F_1}{F_0}$.

\medskip

Observe that $\tau^2(F_0)=F_0$.
This implies that
the $\tau$-degree of
every nonzero monomial (with respect to variables $x_0,\dots,x_{n-1})$
of $F_0$ is divisible by $3$.
This means that in the $\tau$-decomposition of $F_0$ there are only
components with $\tau$-degrees $0$ and $3$.
Let
$
F_0=v_0+v_3,
$
where $v_0\in k[X]$ is  $\tau$-homogeneous with
$\deg_{\tau}(v_0)=0$ (that is, $\tau(v_0)=v_0$),
and $v_3\in k[X]$ is  $\tau$-homogeneous
with $\deg_{\tau}(v_3)=3$ (that is, $\tau(v_3)=\E^3(v_3)=-v_3$).
Of course $d(v_0)=d(v_3)=0$.
Observe  that $\varrho(F_0)=F_0$.
Hence,
$$
v_0+v_3=F_0=\varrho(F_0)=\varrho(v_0)+\varrho(v_3).
$$
Since the $\tau$-decomposition of $F_0$ is unique, we deduce
(by Proposition \ref{TAUrho}), that
$v_3=\varrho(v_0)$ and $v_0=\varrho(v_3)$,
and so,
the $\tau$-decomposition of $F_0$ is of the form
$
F_0=v_0+\varrho(v_0).
$
Moreover,
$F_1=\tau(F_0)=\tau(v_0)+\tau(\varrho(v_0))=v_0+\E^3\varrho(v_0)=v_0-\varrho(v_0)$.

\smallskip

We do a similar procedure with the polynomial $G_0$.
We  first observe that $\tau^3(G_0)=G_0$, and $\varrho(G_0)=-G_0$, and then
we obtain the following three $\tau$-decompositions:
$
G_0=r_0-\varrho(r_0)+\varrho^2(r_0),\quad
G_1=r_0-\E^2\varrho(r_0)+\E^4\varrho^2(r_0),\quad
G_2=r_0-\E^4\varrho(r_0)+\E^2\varrho^2(r_0),
$
where $r_0$ is homogeneous polynomial of degree $2$ which is $\tau$-homogeneous
of $\tau$-degree zero.
Consider now the rational functions $g_1,g_2,g_3\in k(X)$ defined by
$$
\textstyle
g_1=\frac{\varrho(v_0)}{v_0},\quad
g_2=\frac{\varrho(r_0)}{r_0},\quad
g_3=\frac{\varrho^2(r_0)}{r_0}.
$$
These functions are $\tau$-homogeneous.
They are homogeneous of degree zero (in the ordinary sense)
and they are constants of $d$. Moreover,
the quotients $\frac{w_1}{w_0}$, \ $\frac{w_2}{w_0}$, \ $\frac{w_3}{w_0}$,
belong to $k(g_1,g_2,g_3)$.
In fact:
$$
\begin{array}{lcl}
\frac{w_1}{w_0}&=&\frac{F_1G_1}{F_0G_2}
=\frac{\left(v_0-\varrho(v_0)\right)
\left(r_0-\E^2\varrho(r_0)+\E^4\varrho^2(r_0)\right)}
{\left(v_0+\varrho(v_0)\right)
\left(r_0-\E^4\varrho(r_0)+\E^2\varrho^2(r_0)\right)}
=\frac{v_0^{-1}r_0^{-1}\left(v_0-\varrho(v_0)\right)
\left(r_0-\E^2\varrho(r_0)+\E^4\varrho^2(r_0)\right)}
{v_0^{-1}r_0^{-1}\left(v_0+\varrho(v_0)\right)
\left(r_0-\E^4\varrho(r_0)+\E^2\varrho^2(r_0)\right)}\medskip\\
&=&
\frac{(1-g_1)(1-\E^2g_2+\E^4g_3)}
{(1+g_1)(1-\E^4 g_2+\E^2 g_3)},
\end{array}
$$
and so, $\frac{w_1}{w_0}\in k(g_1,g_2,g_3)$.
By a similar way we show that
$\frac{w_2}{w_0}$ and $\frac{w_3}{w_0}$ also belong to
$k(g_1,g_2,g_3)$.
Hence, by Proposition \ref{HomDerZero}, the elements $g_1,g_2,g_3$ are algebraically
independent over $k$ and $k(X)^{E,d}=k(g_1,g_2,g_3)$.
It follows from Proposition \ref{TAUpropGlowne}, that
for each $g_j$ there exists
a homogeneous rational function $f_j\in k(Y)$ such that
$\Delta(f_j)=0$ and $@(f_j)=g_j$.
We know, by Theorem \ref{RCtw}, that the elements
$v, f_1,f_2,f_3$, are algebraically independent over $k$,
and $k(Y)^\Delta=k(v,f_1,f_2,f_3)$.
\end{dowod}

Now we assume that $p>q$ are primes, and $n=pq$.
In the above proof we used the  explicit  form of the cyclotomic polynomial
$\Phi_6(t)$.
Let $\Phi_{pq}=\sum c_jt^j$.
In $1883$, Migotti \cite{Mig} showed that all $c_j$ belong to $\{-1,0,1\}$.
In $1964$ Beiter \cite{Bei} gave a criterion on $j$ for $c_j$ to be $0$, $1$ or $-1$.
A similar result, but more elementary, gave in $1996$,
Lam and Leung \cite{Lam}. Their criterion is based on the fact that
$\f(pq)=(p-1)(q-1)$
can be expressed uniquely in the form
$rp+sq$
where $r,s$ are nonnegative integers. Thus, we have the equality
$$
\f(pq)=rp+sq\quad\mbox{with}\quad r,s\in\mathbb N.
$$
The numbers $r,s$ are uniquely determined, and it is clear that
$0\leqslant r\leqslant q-2$, \ $0\leqslant s\leqslant p-2$, \  $r=r_1-1$ and $s=s_1-1$,
where $r_1\in\{1,\dots,q-1\}$, $s_1\in\{1,\dots,p-1\}$ such that $\MOD{r_1p}{1}{q}$
and $\MOD{s_1q}{1}{p}$.
Using the numbers $r,s$, Lam and Leung proved:

\bblem{\cite{Lam}}{PQlem}
Let $\Phi_{pq}(t)=\sum_{k=0}^{\f(pq)} c_k t^k$. Then
$$
\begin{array}{lcl}
c_k=1&\iff& k=ip+jq, \ i\in\{0,1,\dots,r\}, \ j\in\{0,1,\dots,s\}; \smallskip\\
c_k=-1&\iff& k=ip+jq+1, \ i\in\{0,1,\dots,(q-2)-r\}, \ j\in\{0,1,\dots,(p-2)-s\}.
\end{array}
$$
\elem

Now we may prove the following theorem.

\btw{PQtw}
If $n=pq$ where $p>q$ are primes, then
$$
k(Y)^\Delta=k(v, f_1,\dots,f_{m-1})
$$
with $m=p+q-1$, where $v=y_0\cdots y_{n-1}$ and $f_1,\dots,f_{m-1}\in k(Y)$
are homogeneous rational functions such that
$v,f_1,\dots,f_{m-1}$ are algebraically independent over $k$.
\etw

\begin{dowod}
We use the same idea as in the
proofs of Theorem \ref{RPStw} and Example \ref{PQexa}.
We have: $\f(n)=(p-1)(q-1)$ and $m=n-\f(n)=p+q-1$.
For each $i\in\mathbb Z$, let us denote:
$$
F_i=\prod_{j=0}^{p-1} u_{jq+i},\quad G_i=\prod_{j=0}^{q-1} u_{jp+i}.
$$
In particular, $F_0=u_0u_qu_{2q}\cdots u_{(p-1)q}$
 \ $G_0=u_0u_pu_{2p}\cdots u_{(q-1)p}$.
Observe that if $i=bq+c$, where $b,c\in\mathbb Z$ and $0\leqslant c<q$, then $F_i=F_c$.
Similarly, if $i=bp+c$, where $b,c\in\mathbb Z$ and $0\leqslant c<p$, then $G_i=G_c$.
Let $A$ be the set of all indexes $k\in\{0,1,\dots,\f(pq)\}$ with $c_k=1$,
and let $B$ be the set of all indexes $k\in\{0,1,\dots,\f(pq)\}$ with $c_k=-1$.
It is clear that $A\cap B=\emptyset$, $A\neq\emptyset$, $B\neq\emptyset$, and
$w_0=\frac{N}{D}$ where
$N=\prod_{k\in A}u_k$, \ $D=\prod_{k\in B}u_k$.
It follows from Lemma \ref{PQlem}, that
$$
N=\prod_{i=0}^r\prod_{j=0}^s u_{ip+jq},\quad
D=\prod_{i=0}^{(q-2)-r} \ \prod_{j=0}^{(p-2)-s} u_{ip+jq+1}.
$$
It is easy to check that $\prod_{i=0}^r F_{ip}=N\cdot S$
 \ and \ $\prod_{j=0}^{p-2-s} G_{jq+1}=D\cdot T$, where
$$
S=\prod_{i=0}^r\prod_{j=s+1}^{p-1} u_{ip+jq}
\quad\mbox{and}\quad
T=\prod_{j=0}^{p-2-s}\prod_{i=q-2=r+1}^{q-1} u_{ip+jq+1}
$$
Now we will show that $S=T$. First observe that $S$ and $T$ have the
same number of factors, which is equal to $(r+1)(p-s-1)$.
Next observe that
$$
S=\prod_{i=0}^r\prod_{j=0}^{p-s-2} u_{ip+(s+1+j)q}
\quad\mbox{and}\quad
T=\prod_{j=0}^{p-2-s}\prod_{i=0}^{r} u_{(q-r-1+i)p+jq+1}.
$$
Thus, it is enough to show that, that for $i\in\{0,\dots,r\}$
and $j\in\{0,1,\dots,p-s-2\}$,
we have
$
\MOD{(s+1+j)q+ip}{(q-r-1+i)p+jq+1}{pq}.
$
But it is obvious, because $(p-1)(q-1)=rp+sq$. Therefore, $S=T$ and we have
$$
w_0=\frac{\prod_{i=0}^r F_{ip}}{\prod_{j=0}^{p-2-s}G_{jq+1}}.
\leqno{(\ast)}
$$

\medskip

Now we do exactly the same as in the proof of Example \ref{PQexa}.
We have  the homogeneous polynomials $F_0,\dots,F_{q-1}$ and $G_0,\dots,G_{p-1}$,
which are constants of $d$, and $F_i=\tau^i(F_0)$, $G_i=\tau^i(G_0)$,
$\deg F_i=p$, $\deg G_i=q$, for each $i$.
Observe that $\tau^q(F_0)=F_0$.
This implies that
the $\tau$-degree of
every nonzero monomial (with respect to variables $x_0,\dots,x_{n-1})$
of $F_0$ is divisible by $p$.
This means that in the $\tau$-decomposition of $F_0$ there are only
components with $\tau$-degrees $0,p,2p,\dots,(q-1)p$.
Let
$F_0=\sum_{i=0}^{q-1} v_i$,
where each $v_i$ is a $\tau$-homogeneous polynomial from $k[X]$, and
$\tau(v_i)=\E^{pi}v_j$.
Of course $d(v_i)=0$ for all $i$ (because $\tau d=\E d\tau$),
and  $\deg (v_i)=p$.
But  $\varrho(u_j)=\E^{-j}u_j$ (see Lemma \ref{NVlemC}), so
$\varrho(F_0)=\pm F_0$, and we have
$$
v_0+v_1+\dots+v_{m-1}=F_0=\pm \varrho(F_0)=\pm(\varrho(v_0)\pm\varrho(v_1)\pm
\dots\pm\varrho(v_{m-1})
$$
Since the $\tau$-decomposition of $F_0$ is unique,
we deduce (by Proposition \ref{TAUrho}), that
$
v_1=\pm\varrho(v_0),\quad\ v_2=\pm\varrho(v_1),\quad \dots,\quad  v_{m-1}=\pm\varrho(v_{m-2}),\quad v_0=\pm\varrho(v_{m-1}),
$
and we have  $v_j=\pm \varrho^j(v_0)$ for all $j=0,1,\dots,q-1$.
Therefore, the $\tau$-decomposition of $F_0$ is of the form
$
F_0=v_0+\sum_{i=1}^{q-1} b_i\varrho^i(v_0),
$
where $b_1,\dots,b_{m-1}\in\{-1,1\}$. This implies that
$
F_1=\tau(F_0)=v_0+\sum b_i\E^{ip}\varrho(v_0).
$
We do the same for $F_2=\tau(F_1)=\tau^2(F_0)$, and for all $F_j$.
Thus,
for all $j=0,1,\dots,m-1$, we have
$$
F_j=v_0+\sum_{i=1}^{q-1}c_{ji}\varrho^i(v_0),
$$
where each $c_{ji}$ belongs to the ring $\mathbb Z[\E]$.
We do a similar procedure with the polynomial $G_0$.
First  observe that $\tau^p(G_0)=G_0$
and $\varrho(G_0)=\pm G_0$, and then
we obtain $\tau$-decompositions of the forms
$$
G_j=r_0+\sum_{i=1}^{p-1} b_{ji}\varrho^i(r_0),
$$
where each $c_{ji}$ belongs to $\mathbb Z[\E]$.
where $r_0$ is a homogeneous polynomial of degree $q$
which is $\tau$-homogeneous of $\tau$-degree zero.

Consider now the elements $g_1,\dots,g_{m-1}\in k(X)$
defined by
$$
g_i=\frac{\varrho^i(v_0)}{v_0},\quad
g_{q-1+j}=\frac{\varrho^j(r_0)}{r_0},
$$
for $i=1,\dots,q-1$, and $j=1,\dots,p-1$.
These elements are $\tau$-homogeneous.
They are homogeneous of degree zero (in the ordinary sense)
and they are constants of $d$.
We know, by the above construction, that each element of the form
$
\frac{1}{v_0}\tau^i(F_j)\quad\mbox{or}\quad
\frac{1}{r_0}\tau^i(G_j)
$
belongs to the field $k(g_1,\dots,g_{m-1})$.
But, by $(\ast)$, for each $a=0,\dots,m-1$, we have
$$
w_a\frac{r_0^{p-1-s}}{v_0^{r+1}}=
\tau^a(w_0)\frac{r_0^{p-1-s}}{v_0^{r+1}}=
\frac{\prod_{i=0}^r \frac{\tau^a(F_{ip})}{v_0}}
{\prod_{j=0}^{p-2-s}\frac{\tau^a(G_{jq+1})}{r_0}},
$$
and hence, each element
$w_a r_0^{p-1-s} v_0^{-(r+1)}$ belongs to $k(g_1,\dots,g_{m-1})$.
This implies, that for every $j-1,\dots,m-1$, the quotient
$$
\frac{w_j}{w_0}=\frac{r_0^{p-1-s} v_0^{-(r+1)}w_j}
{r_0^{p-1-s} v_0^{-(r+1)} w_0}
$$
belongs to $k(g_1,\dots,g_{m-1})$.
Hence, by Proposition \ref{HomDerZero}, the elements $g_1,\dots,g_m$
are algebraically
independent over $k$ and $k(X)^{E,d}=k(g_1,\dots,g_{m-1})$.
It follows from Proposition \ref{TAUpropGlowne}, that
for each $g_j$ there exists
a homogeneous rational function $f_j\in k(Y)$ such that
$\Delta(f_j)=0$ and $@(f_j)=g_j$.
We know, by Theorem \ref{RCtw}, that the elements
$v, f_1,\dots,f_{m-1}$, are algebraically independent over $k$,
and $k(Y)^\Delta=k(v,f_1,\dots,f_{m-1})$.
This completes our proof of Theorem \ref{PQtw}.
\end{dowod}

We already know a structure of the field $k(Y)^\Delta$ but only in the following
two cases,
when $n$ is a power of a prime number (Theorem \ref{RPStw}),
and when $n$ is the product of two
prime numbers (Theorem \ref{PQtw}). We do not know what happens in all other cases.
Is this field always a purely transcendental extension of $k$\,?
What is in the cases $n=12$ or $n=30$ or $n=105$\,?


\end{document}